\documentclass[MSNbibl,number,citesort,seceqns,dvips]{arxbj}

\aid{0}
\volume{18}
\issue{3}
\pubyear{2012}
\firstpage{945}
\lastpage{974}
\doi{10.3150/11-BEJ364}

\makeatletter
\newcommand{\eqref}[1]{(\ref{#1})}
\newtheorem{theorem}{Theorem}[section]
\newtheorem{lemma}[theorem]{Lemma}
\newproclaim{definition}[theorem]{Definition}
\newremark{remark}{Remark}

\makeatother

\begin{document}
\begin{frontmatter}

\title{The log-linear group-lasso estimator and its asymptotic properties}
\runtitle{The log-linear group-lasso estimator}

\begin{aug}
%%%% inicialai - be tarpu
\author[a]{\fnms{Yuval} \snm{Nardi}\thanksref{a}\ead[label=e1]{ynardi@ie.technion.ac.il}}
\and
\author[b]{\fnms{Alessandro} \snm{Rinaldo}\corref{}\thanksref{b}\ead[label=e2]{arinaldo@stat.cmu.edu}}
\runauthor{Y. Nardi and A. Rinaldo}
\address[a]{Faculty of Industrial Engineering and Management,
Technion--Israel Institute of Technology,
Haifa 32000, Israel.
\printead{e1}}

\address[b]{Department of Statistics,
Carnegie Mellon University,
Pittsburgh, PA 15213-3890, USA.\\
\printead{e2}}
\end{aug}

% HISTORY:
\received{\smonth{10} \syear{2009}}
\revised{\smonth{12} \syear{2010}}

% ABSTRACT
%
\begin{abstract}
We define the group-lasso estimator for the natural parameters of the
exponential families of distributions representing hierarchical
log-linear models under multinomial sampling scheme. Such estimator
arises as the solution of a convex penalized likelihood optimization
problem based on the group-lasso penalty. We illustrate how it is
possible to construct an estimator of the underlying log-linear model
using the blocks of nonzero coefficients recovered by the group-lasso procedure.
We investigate the asymptotic properties of the group-lasso estimator
as a~model selection method in a double-asymptotic framework, in which
both the sample size and the model complexity grow simultaneously.
We provide conditions guaranteeing that the group-lasso estimator is
model selection consistent, in the sense that, with overwhelming
probability as the sample size increases, it correctly identifies all
the sets of nonzero interactions among the variables.
Provided the sequences of true underlying models is sparse enough,
recovery is possible even if the number of cells grows larger than the
sample size. Finally, we derive some central limit type of results for
the log-linear group-lasso estimator.
\end{abstract}

% KEYWORDS
%

\end{frontmatter}

%s1 ###
\section{Introduction}\label{sec1}

The theory of log-linear models has produced a variety of statistical
methodologies and theoretical results for the analysis of categorical
data that have found applications in numerous scientific areas, ranging
from social and biological sciences, to medicine, disclosure limitation
problems, data-mining, image analysis, finger-printing, language
processing and genetics.

%Statistical analysis of categorical data via log-linear models
%Log-linear model analysis of categorical data relies on an important
%and widespread set of statistical methodologies that have found
%applications in very diverse scientific areas, ranging from social
%and biological sciences, to medicine, disclosure limitation problems,
%data-mining, image analysis, finger-printing, language processing and
%genetics.
Inherently, log-linear modeling is a model selection procedure for
contingency tables that encompasses testing a number of statistical
models for the joint distribution of a set of categorical variables.
The classical asymptotic theory of model selection and goodness-of-fit
testing is well developed and understood for the `small $p$ and large
$N$' case, that is, the case in which the sample size $N$ is much
larger than the number $p$ of candidate parameters. It is applicable to
a variety of goodness-of-fit measures, such as Pearson's $\chi^2$, the
likelihood ratio statistic and, more generally, any statistics
belonging to the power-divergence family of Read and Cressie \cite{CR88}.
The applicability and validity of these methods demand the availability
of large sample sizes and the existence of the maximum likelihood
estimate (MLE).

In recent years, the importance and usage of log-linear modeling
methodologies have increased dramatically with the compilation and
diffusion of large databases in the form of sparse contingency tables.
In such instances, the number of sampled units is not much different,
in fact often smaller, than the number of cells, so that most of the
cell entries are very small or zero counts. In high-dimensional
settings, the traditional methodologies indicated above are inadequate.
First off, the number of log-linear models grow extremely fast with the
number of variables (for example, there are 7580 hierarchical models
for a 5-way table!), and selecting an optimal model involves exploring
a space of models of virtually infinite dimension. Secondly, for a
given model of even moderate complexity, under a sparse scenario, the
MLE is unlikely to exist. This implies that the information content
present in the data is not sufficient to estimate all the parameters of
the model, and, therefore, the possibility for inference is only
limited to portions of the parameter space (see Rinaldo, Fienberg and Zhou \cite{ERGM} for
details). As a result, traditional goodness-of-fit testing and
model selection will produce very poor, if not completely erroneous,
asymptotic approximations. It is quite clear that a more appropriate
statistical formalization requires the consideration of a `large $p$' setting.

In this article, we study a methodology for log-linear model selection
that is particularly suited to high-dimensional tables, and we describe
some of its asymptotic properties. Our results are akin to the
asymptotic optimality of the lasso estimator in high dimensional least
squares problems, where the recovery of the sparsity pattern of an
unknown set of parameters in noisy settings via $\ell_1$-regularization
is possible, even if the number of parameters grows faster than the
sample size. See, in particular, Meinshausen and B\"{u}hlmann \cite{BM06}, Zhao and Yu \cite{ZY06}, Wainwright~\cite
{W06} and, for a different approach, Greenshtein \cite{GR06} and Greenshtein and Ritov \cite{GRR04}.
Existing work on penalized likelihood problems involving $\ell
_1$-regularization for discrete data include the nonasymptotic results
about $\ell_2$ consistency for estimation in high-dimensional
generalized linear models via the lasso by van de Geer \cite{sara106,sara206},
and the analysis by Wainwright, Ravikumar and Lafferty~\cite{WPJ06} on the consistency of $\ell
_1$-regularized logistic regression with binary variables under a
double asymptotic framework.
In Section \ref{secconclusions}, we discuss in detail the differences
between our problem and solutions and the existing results.

We formulate the log-linear model selection problem as a convex
penalized maximum likelihood problem based on the group-lasso, a convex
penalty function introduced by Yuan and Lin~\cite{YL06} in a nonasymptotic ANOVA
setting and further analyzed by Nardi and Rinaldo \cite{NR08}. The group-lasso
regularization is an extension of the lasso, or $\ell_1$, penalty
function designed to penalized groups of coefficients simultaneously.
It has been shown to be effective in logistic regression problems by
Meier, van der Geer, and B\"{u}hlmann \cite{MVB06} and has been used in applications involving log-linear
modeling of sparse contingency tables in Dahinden et al. \cite{LOGLASSO06}.

The paper is organized as follows. In Section \ref{secloglin}, we
describe the log-linear model settings we will be considering. The
direct sum decomposition of the natural parameter space by log-linear
subspaces defines a partition of the parameters in blocks of different
dimensions, which are utilized as arguments of the group penalty
function. In Section~\ref{secgrouplasso}, we describe the group-lasso
estimator for log-linear models, which can be computed by solving a~convex program. %and can be interpreted as a smoothed MLE (see Section
%of log-linear models,
Next, we show that the group-lasso estimator produces, in turn, an
estimator of the underlying log-linear model, which is constructed
simply by isolating the nonzero blocks of the group-lasso estimates.
Section \ref{secasymptotics} outlines our contribution by studying the
consistency properties of the group-lasso estimator as a model
selection procedure. We formulate a general double-asymptotic framework
in which we allow both the sample size and the model complexity to
grow. In Section \ref{secmodselconsistency}, we derive conditions
guaranteeing that the model estimates are consistent, that is,
asymptotically, the group-lasso correctly identifies the set of
interactions making up the underlying model. We conclude our analysis
with some central limit results in Section \ref{secclt}.
The proofs appear in Section~\ref{secproofs}.\looseness=-1 %and the appendix
%contain additional material on the linear subspaces used for
%hierarchical log-linear models.

%s1.1 ###
\subsection{Notation}
Let $X_1, \ldots, X_K$ denote $K$ categorical variables, where each
$X_k$ takes values in $\mathcal{I}_k=\{1, \ldots, I_k\}$, with
$I_k\geq
2$ an integer. Set $\mathcal{I} = \mathcal{I}_1 \times\cdots\times
\mathcal{I}_K$. We denote by $\mathbb{R}^{\mathcal{I}}$ the class of
real-valued functions on $\mathcal{I}$ which is the vector space of
real-valued $K$-dimensional arrays indexed by the multi-index~$\mathcal
{I}$. The vector space $\mathbb{R}^{\mathcal{I}}$ can be naturally
represented as a Euclidean space of dimension $I \equiv\prod_{k=1}^K
I_k$. This identification can be realized in a~straightforward fashion
by ordering $\mathcal{I}$ as a linear list using any bi-jection between
$\mathcal{I}$ and $\{1,2,\ldots,I\}$.
Each element $i$ of $\mathcal{I}$, called a \textit{cell}, is a
multi-index $i = (i_1, \ldots, i_K)$. %, whose $k$-th coordinate is the
%valued taken on by the $k$-th variable .
%Therefore, each cell $i$ can be though of either as a multi-index or
%as a number between $1$ and $I$.
Using this coordinate vector representation, for any array $\mathbf{x}
\in
\mathbb{R}^\mathcal{I}$, $\mathbf{x}_i$ is the number indexed by the
coordinate $i \in\mathcal{I}$. Also, the standard inner product
$\langle\mathbf{x}, \mathbf{y} \rangle= \sum_{i \in\mathcal{I}}
\mathbf{x}_i
\mathbf{y}_i$ and the induced Euclidean norm are well defined for all
$\mathbf{x}, \mathbf{y} \in\mathbb{R}^{\mathcal{I}}$.

We use the following notational conventions. Let $\{x_s, s \in\mathcal
{S} \}$ be a set of vectors of possibly different dimensions, indexed
by some finite set $\mathcal{S}$. We will denote by $\operatorname{vec}\{
x_s,\allowbreak s
\in\mathcal{S}\}$
the vector obtained by staking the $x_s$'s one on top of each others in
the same order as the elements of~$\mathcal{S}$.
For any $d$-dimensional Euclidean vector $x$, we write $\exp
(x)=\operatorname{vec}\{\mathrm{e}^{x_i} , i = 1,\ldots, d \}$, and $\log x=\operatorname
{vec}\{\log
x_i , i= 1,\ldots, d\}$. For a linear subspace~$\mathcal{A}$ of
the $d$-dimensional Euclidean space, we denote with $\mathcal{A}^\bot$
the orthogonal complement of $\mathcal{A}$. If~$\mathcal{B}$ is another
linear subspace orthogonal to~$\mathcal{A}$, we write $\mathcal{A}
\oplus\mathcal{B}$ for the linear subspace obtained as their direct sum.
Similarly, for matrices $\mathrm{U}_1,\ldots, \mathrm{U}_n$ with the
same number
of rows $r$ and number of columns $c_1,\ldots,c_n$, respectively, we
denote the operation of adjoining them into one matrix of dimension $r
\times\sum_k c_k$ with $
\bigoplus_{k=1}^n \mathrm{U}_k = [ \mathrm{U}_1 \cdots\mathrm{U}_n ]$.

Throughout the article, we will consider random vectors and functions
of random vectors whose probability distributions will always be clear
from the context. As a result, we will use the generic notation
$\mathbb
{P}(\mathcal{O})$ for the probability of an event $\mathcal{O}$ defined
by such vectors and will write $\mathbb{E}$ for the corresponding
expectation operator.

%s2 ###
\section{Log-linear models}\label{secloglin}

We will be considering the usual log-linear modeling setting, which we
now describe. See (see also Bishop, Fienberg and Holland \cite{BFH75}, Haberman \cite{HAB74},
Lauritzen \cite{LAU96}). The $K$
categorical random variables $X_1, \ldots, X_K$ have an unknown joint
distribution given by the strictly positive probability vector $\bolds
{\pi }$ in $\mathbb{R}^\mathcal{I}$ with coordinates
\[
\bolds{\pi}_{i_1,\ldots,i_K} = \mathbb{P} \bigl( (X_1, \ldots, X_K) =
(i_1, \ldots, i_K)\bigr), \qquad(i_1, \ldots, i_K) \in\mathcal{I}.
\]
The positivity of $\bolds{\pi}$ is a crucial assumption, ruling out the
case of \textit{structural zeros}, that is, cells that can never be observed.

We observe $N$ independent and identically distributed realizations of
the random vector $(X_1, \ldots, X_K)$. Their cross-classification
results in a random integer-valued vector $\mathbf{n} \in\mathbb
{R}^\mathcal{I}$, called a \textit{contingency table}, whose $i$th
coordinate entry $\mathbf{n}_{i}$ corresponds to the number of times the
cell combination $i = (i_1, \ldots, i_K)$ was observed in the sample.
The table $\mathbf{n}$ has a $\operatorname{Multinomial}(N, \bolds{\pi})$
distribution.

Log-linear model theory is concerned with drawing inferences on $\bolds
{\pi}$ based on the observed table $\mathbf{n}$. %i.i.d. sample of
%size
%$N$ for $(X_1, \ldots, X_K)$.
Specifically, let $\mathbf{m} = \mathbb{E} \mathbf{n} = N \bolds
{\pi}$ denote
the (necessarily positive) cell mean vector of $\mathbf{n}$ and set
$\bolds{\mu} = \log\mathbf{m}$. Notice that estimating $\bolds{\mu
}$ is equivalent
to estimating $\bolds{\pi}$. A log-linear model is specified by
prescribing a linear subspace $\mathcal{M}$ of $\mathbb{R}^\mathcal{I}$
containing the constant functions and then requiring that $\bolds{\mu}$
belongs to $\mathcal{M}$. Indeed, any point in $\mathcal{M}$ represents
a different cell mean vector, and, therefore, a different probability
distribution over $\mathcal{I}$.

%We $N$ independent and identically distributed realizations of the
%$K$-dimensional random vector $(X_1, \ldots, X_K)$ from an unknown

%Letting $\mathcal{I} = \bigotimes_{k=1}^K \mathcal{I}_k$, the set $
%the vector space of real-valued $K$-dimensional arrays indexed by $
%multi-index $i = (i_1, \ldots, i_K)$, whose $k$-th coordinate is the
%valued taken on by the $k$-th variable .
%The vector space can be naturally represented as a the Euclidean space
%of dimension $I = \prod_k I_k$. Therefore, each cell $i$ can be though
%of either as a multi-index or as a number between $1$ and $I$.
%Using this notation, for any $\mathbf{x} \in\mathbb{R}^\mathcal{I}$, ${
%standard inner product $\langle\mathbf{x}, \mathbf{y} \rangle= \sum_{i
%${
% Finally, for a vector $\mathbf{x} \in\mathbb{R}^d$, we will write $
%e^{{

%For convenience, we will identify $\mathbb{R}^{\mathcal{I}}$ with the
%Euclidian space $\mathbb{R}^{I}$, where $I = \prod_k I_k$. This
%identification can be made in a straightforward fashion by ordering $

%In log-linear model theory, the joint distribution of $(X_1, \ldots,
%X_K)$ is fully specified by representing the cell mean vector $
%= \mathbb{E} \mathbf{n} = N \mbf{\pi}$ by means of certain linear
%subspace $\mathcal{M}$ of $\mathbb{R}^I$ containing $\log\mathbf{m}$.
%%,
%to the extent that log-linear models themselves are defined by such
%subspaces.
%Namely, by fixing $\mathcal{M}$, it follows that the logarithms of the
%cell mean vectors must satisfy specific linear constraints, to be
%specified below, which completely characterize the underlying
%distribution.
Then, for a given table $\mathbf{n}$ the log-likelihood function $\ell^*$
at $\bolds{\mu} \in\mathcal{M}$ is (see Haberman~\cite{HAB74}, page~11)
\[
\ell^*(\bolds{\mu}) = \cases{
\displaystyle\sum_{i \in\mathcal{I}} \mathbf{n}_i \log\frac{\mathbf
{m}_i}{\langle
\mathbf{m}, \mathbf{1} \rangle} + \log N! - \displaystyle\sum_{i \in\mathcal
{I}} \log
\mathbf{n}_i!, & \quad$\mbox{ if } \langle\mathbf{m}, \mathbf{1} \rangle
= N,$\vspace*{2pt}\cr
\mbox{undefined}, & \quad $\mbox{otherwise},$}
\]
where $\mathbf{m} = \exp(\bolds{\mu})$ and $\mathbf{1} \in\mathbb
{R}^{\mathcal
{I}}$ is the $I$-dimensional vector containing ones. Indeed, because of
the Multinomial sampling assumption, $\ell^*$ is only defined over the
nonconvex set $\widetilde{\mathcal{M}} \subsetneq\mathcal{M}$ given by
\[
\widetilde{\mathcal{M}} = \{ \bolds{\mu} \in\mathcal{M} \dvt
\langle\mathbf{m}, \mathbf{1} \rangle= N\}.
\]
This parametrization is clearly quite inconvenient. Fortunately, it is
possible to reparametrize the log-likelihood function as concave
function defined over the entire~$\mathbb{R}^k$, where $k$ is the
dimension of $\widetilde{\mathcal{M}}$. Specifically, let~$ \mathcal
{R}(\mathbf{1})$ be the one-dimensional subspace of $\mathbb
{R}^\mathcal
{I}$ spanned by $\mathbf{1}$ and consider the linear subspace
$\mathcal{M}
\cap\mathcal{R}(\mathbf{1})^{\bot} \subset\mathbb{R}^{\mathcal
{I}}$ of
dimension $k = \operatorname{dim}(\mathcal{M})-1$.

\begin{lemma}\label{lemlik}
Let $\mathrm{U}$ be any full-rank matrix whose columns span $\mathcal{M}
\cap\mathcal{R}(\mathbf{1})^{\bot}$ and consider the function
%
%
%e1 ###
\begin{equation}\label{eqexpfam}
\ell(\theta) = \langle\mathrm{U}^{\top} \mathbf{n}, \theta
\rangle- N \log
\langle\exp(\mathrm{U} \theta), \mathbf{1} \rangle+ \log N! -
\sum_{i \in
\mathcal{I}} \log\mathbf{n}_i!,\quad \theta\in\mathbb{R}^k.
\end{equation}
Then, for each $\widetilde{\bolds{\mu}} \in\widetilde{\mathcal{M}}$
there exists one $\theta\in\mathbb{R}^k$ such that
%
%
%e2 ###
\begin{equation}\label{eqhomeo}
\exp(\widetilde{\bolds{\mu}}) =\frac{N}{\langle\exp(\mathrm{U}
\theta) ,
\mathbf{1} \rangle} \exp(\mathrm{U}\bolds{\theta})\quad \mbox{and}\quad
\ell
(\theta) = \ell^*(\widetilde{\bolds{\mu}})\qquad \mbox{for each }
\mathbf{n},
\end{equation}
and, conversely, for each $\theta\in\mathbb{R}^k$ there exists one
$\widetilde{\bolds{\mu}} \in\widetilde{\mathcal{M}}$ satisfying the
above identities.
\end{lemma}
This reparametrization is essentially equivalent to reduction to
minimal form of the underlying exponential family of distributions for
the cell counts via sufficiency. In fact, the previous display shows
that each log-linear model $\mathcal{M}$ specifies a full, regular
exponential family of dimension $k = \operatorname{dim}(\mathcal{M}) - 1$, natural
sufficient statistic $\mathrm{U}^{\top} \mathbf{n}$ and natural
parameter space
$\mathbb{R}^k$ (see, e.g., Brown~\cite{BROWN89}).
Throughout the article, we take \eqref{eqexpfam} as the operative
definition of log-likelihood function. Notice that the log-likelihood
function depends on the choice of the design matrix $\mathrm{U}$. %See
%the
%next section more details on how to compute such design matrices.

%Finally, we remark that the gradient and Hessian matrix for $\ell(
%one can see that
%( \mathbf{n} - \mathbf{m} )
%and
%where $\mathbf{m} = \frac{N}{\langle\mathbf{b},\mathbf{1}\rangle}
%matrix
%with diagonal $\mathbf{m}$. It is worth pointing out that, because
%these
%models are exponential families, the negative Hessian is the
%covariance matrix of the natural sufficient statistics $\m{U}^{\top} {

%s2.1 ###
\subsection{Log-linear subspaces for hierarchical log-linear
models}\label{secloglinsubs}
Although log-linear models are defined by generic linear manifolds of
$\mathbb{R}^{\mathcal{I}}$, in practice it is customary to consider
only very specific classes of linear subspaces, which are also
characteristic of ANOVA models and experimental design, yielding
hierarchical log-linear. In this section, we briefly describe such
subspaces. See Darroch, Lauritzen and Speed~\cite{DLS80}, Appendix B in Lauritzen~\cite{LAU96} and Rinaldo \cite
{ALEb06} for details.

%Although log-linear models are defined by generic linear manifolds of $
%only very specific classes of linear subspaces, which are also
%characteristic of ANOVA models and experimental design. %These
%subspaces present considerable advantages in terms of interpretability
%and ease of computation and can be constructed easily by exploiting
%various correspondences between combinatorial structure of the power
%set of $\mathcal{K} = \{ 1, \ldots, K\}$ and a certain direct sum
%decomposition of $\mathbb{R}^{\mathcal{I}}$, to be described below.
A rather intuitive way of specifying a certain dependence structure
among the $K$ variables of interest is to provide a list of the
interactions among them.
Then, the associated statistical model is representable as a class of
subsets of $\mathcal{K} \equiv\{ 1,2, \ldots, K \}$, each one
indicating a different type of interaction. In fact, every subset~$h$
of~$\mathcal{K}$ can be given a~straightforward ANOVA-type of an
interpretation, based on its cardinality $|h|$, so that~$h$ identifies
an interaction of order $|h|-1$ among the variables $\{ i \dvt i \in
h \}$. For example, if $|h|=1$, then $h$ is a main effect, if $h =
\varnothing$, then $h$ is the grand mean, and so on.

Formally, let $2^{\mathcal{K}}$ be the power set of $\mathcal{K}$,
which we view as a lattice with respect to the partial order induced by
the operation of taking subset inclusion.

\begin{definition}\label{defhloglin}
A hierarchical log-linear model $\Delta$ is a collection of subsets of
$2^{\mathcal{K}}$ such that $h \in\Delta$ and $h' \subset h$ implies
$h' \in\Delta$. An interaction model $\mathcal{H}$ is just a subset of
$2^{\mathcal{K}}$.
\end{definition}
By definition, once an interaction term is part of $\Delta$, all lower
order interactions are included. %, i.e. the model is hierarchical.
%%Then, hierarchical log-linear model can be defined in a purely
%combinatoric form as a subset of $2^\mathcal{K}$.
Notice that Definition \ref{defhloglin} includes as special case the
class of graphical and hierarchical models (see, e.g., Lauritzen \cite{LAU96}).
Though our analysis is valid also for the larger
class of interaction models, we focus only on hierarchical log-linear
models, primarily because the interpretability of interaction
log-linear models is very limited.

%More generally, a (non-necessarily hierarchical) log-linear model is
%simply a subset $\mathcal{H}$ of $2^\mathcal{K}$. Despite our analysis
%is valid also for the larger class of log-linear models, the
%interpretability of these unrestricted models is unclear.

To any given hierarchical model $\Delta$, there corresponds one
log-linear subspace \mbox{$\mathcal{M}_{\Delta} \subset\mathbb
{R}^{\mathcal
{I}}$}, constructed as the direct sums of subspaces of $\mathbb
{R}^\mathcal{I}$ indexed by the subsets of $\mathcal{K}$ belonging to
$\Delta$. Specifically,
%
%
%e3 ###
\begin{equation}\label{eqMD}
\mathcal{M}_{\Delta} = \bigoplus_{ h \in\Delta} \mathcal{U}_h ,
\end{equation}
where $ \{ \mathcal{U}_h, h \in2^{\mathcal{K}} \}$ are mutually
orthogonal subspaces, called the \textit{subspaces of interactions}. We
refer the reader to Lauritzen \cite{LAU96}, Appendix B, for details on these
subspaces. In particular, $\mathcal{U}_{\varnothing}$ is the
one-dimensional subspace $\mathcal{R}(\mathbf{1})$ and
\[
\operatorname{dim}(\mathcal{U}_h) \equiv d_h = \prod_{k \in h} ( I_k -
1),\qquad h \subseteq\mathcal{K}, h \neq\varnothing.
\]
A design matrix for $\Delta$ can be constructed as follows.
%Below, we describe a class of matrices $\m{U}_h$ whose columns span
%the subspaces $\mathcal{U}_h$, with $h \in2^\mathcal{K}$.
% These bases correspond to contrasts in models of analysis of variance
%and, using Birch's notation \citep[see, in particular, ][]{BFH:75},
%the design matrix for $\mathcal{U}_h$ will encode to the
%the factors in $h$.
For each term $h \subseteq\mathcal{K}$ and factor $k \in\mathcal{K}$,
define the matrix
\[
\mathrm{U}^h_k = \cases{
\mathrm{Z}_k, & \quad$\mbox{if }  k \in h, $\vspace*{2pt}\cr
\mathbf{1}_k, & \quad$\mbox{if }  k \notin h ,$}
\]
where $\mathrm{Z}_k $ is a $I_k \times(I_k-1)$ matrix with entries
%
%
%e4 ###
\begin{equation}\label{eqzk}
\mathrm{Z}_k = \left(
\matrix{
1 & 0 & 0 & 0 & \cdots& 0\vspace*{1pt}\cr
-1 & 1 & 0 & 0 & \cdots& 0\vspace*{1pt}\cr
0 & -1 & 1 & 0 & \cdots& 0\vspace*{1pt}\cr
\vdots& \vdots& \vdots& \vdots& \vdots& \vdots\vspace*{1pt}\cr
0 & 0 & 0 & 0 & \cdots& 1\vspace*{1pt}\cr
0 & 0 & 0 & 0 & \cdots& -1}\right),
\end{equation}
and $\mathbf{1}_k$ is the $I_k$-dimensional column vector of $1$'s. Let
%
%
%e5 ###
\begin{equation}\label{eqUh}
\mathrm{U}_{h} = \bigotimes_{k=1}^K \mathrm{U}^h_k,
\end{equation}
where $\otimes$ denotes the Kronecker product. Then, it is possible to
show that, for each $h \in2^{\mathcal{K}}$, $\mathrm{U}_h$ is a $(I
\times
d_h)$-dimensional full-rank matrix whose columns span $\mathcal{H}_h$.
See Rinaldo \cite{ALEb06}, Section 3, for details.
Thus, the columns of
%
%
%e6 ###
\begin{equation}\label{eqUD}
\mathrm{U}_{\Delta} = \bigoplus_{h \in\Delta, h \neq\varnothing}
\mathrm{U}_h%
\end{equation}
span $\mathcal{M}_\Delta\cap\mathcal{R}(\mathbf{1})^\bot$, and,
therefore, $\mathrm{U}_\Delta$ is a full-rank design matrix for the
log-linear model $\Delta$.
By the same token, the columns of the matrix
\[
\mathrm{U} = \bigoplus_{h \in2^{\mathcal{K}}, h \neq\varnothing}
\mathrm{U}_h
\]
span the $(I-1)$-dimensional subspace $\mathbb{R}^{\mathcal{I}}\cap
\mathcal{R}^\bot$. As a result, any point $\bolds{\mu} \in\mathbb
{R}^{\mathcal{I}}\cap\mathcal{R}(\mathbf{1})^\bot$ can be written as
\[
\bolds{\mu} = \mathrm{U} \theta= \sum_{h \in2^{\mathcal{K}}, h
\neq\varnothing
} \mathrm{U}_h \theta_h,
\]
for some vector
%
%
%e7 ###
\begin{equation}\label{eqthetavec}
\theta= \operatorname{vec}\{ \theta_h, h \in2^{\mathcal{K}}, h \neq
\varnothing\} \in\mathbb{R}^{I-1},
\end{equation}
where $\theta_h$ denotes the $d_h$-dimensional sub-vector of $\theta$
corresponding to the sub-matrix~$\mathrm{U}_h$.\vadjust{\goodbreak}
%and, by Lemma \ref{propdimm}, $d_h = \prod_{k \in h} ( I_k - 1
%)$.

\begin{remark*}
Throughout the document, we will be assuming that the elements of
$2^{\mathcal{K}}$ are ordered in some predefined way, and that any
indexing by subsets of $\mathcal{K}$ is done accordingly. Then, using
such ordering, any vector $\theta\in\mathbb{R}^{I-1}$ can be uniquely
represented like in \eqref{eqthetavec}. Furthermore, we will use the
notation $\sum_h$ to denote summation over all $h \subseteq\mathcal
{K}$ with $h \neq\varnothing$.
\end{remark*}

%s3 ###
\section{The group-lasso estimator for log-linear models}\label{secgrouplasso}
In this section, we define the group-lasso estimator for the set of
interactions of a given hierarchical log-linear model specified by a
subspace $\mathcal{M}_\Delta$.

%Following the results in the previous section, the columns of the
%matrix
%span $(I-1)$-dimensional subspace $\mathbb{R}^{\mathcal{I}}\cap
%and, by Lemma \ref{propdimm}, $d_h = \prod_{k \in h} ( I_k - 1
%)$.
%Accordingly, for any vector $\theta\in\mathbb{R}^{I-1}$, we can write
%where $\theta_h$ denotes the $d_h$-dimensional vector of $\theta$
%corresponding to the sub-matrix $\m{U}_h$.

Using (\ref{eqexpfam}), the log-likelihood function for the saturated
$(I-1)$-dimensional log-linear model is
%
%
%e8 ###
\begin{eqnarray}\label{eqsat}
\ell(\theta) &=& \sum_{h \in2^{\mathcal{K}}, h \neq\varnothing}
\langle\mathrm{U}_h^{\top} \mathbf{n}, \theta_h \rangle- N \log
\biggl\langle\exp\biggl( \sum
_{ h \in2^{\mathcal{K}}, h \neq\varnothing} \mathrm{U}_h \theta_h
\biggr), \mathbf{1}
\biggr\rangle
\nonumber
\\[-8pt]
\\[-8pt]
\nonumber
&&{}+ \log N! - \sum_i \log\mathbf{n}_i!,\qquad \theta\in
\mathbb{R}^{I-1}.
\end{eqnarray}
Notice that the one-dimensional sub-space $\mathcal{R}(\mathbf{1})$
corresponding to the empty set is not included, because of the
multinomial sampling restriction.

For any nontrivial model $\Delta$ with corresponding log-linear
subspace $\mathcal{M}_\Delta$ (i.e., a model different than $\{
\varnothing\}$, which encodes the uniform distribution over $\mathcal
{I}$), let
%
%
%e9 ###
\begin{equation}\label{eqH}
\mathcal{H} = \mathcal{H}(\Delta) =\{ h \dvt h \in\Delta, h \neq
\varnothing\},
\end{equation}
be the collections of sets representing all the interactions in $\Delta
$, or, equivalently, the collections of factor interaction subspaces of
$\mathcal{M}_{\Delta}$, so that $
\operatorname{dim}(\mathcal{M}_\Delta) - 1= \sum_{h \in\mathcal{H}} d_h
\equiv
d_{\mathcal{H}}$. (Notice that~$\mathcal{H}$ differs from $\Delta$ only
because it does not contain the empty set).
We will embed the natural parameter space of $\Delta$, i.e., $\mathbb
{R}^{d_{\mathcal{H}}}$, as a linear subspace of $\mathbb{R}^{I-1}$
consisting of all vectors such that
\[
\cases{
\| \theta_h \| > 0, & \quad$h \in\mathcal{H}$,\vspace*{2pt}\cr    %, \mbox{ with strict
%inequality if } h \mbox{is a maximal subset}\\
\| \theta_h \| = 0, & \quad$h \notin\mathcal{H}.$}
\]
The log-likelihood function for this model is still given by equation
(\ref{eqsat}), where the summations are now taken over the sets $h$ in
the class $\mathcal{H}$.

Let $\Delta_0$ denote the true underlying log-linear model. Thus, there
exists a vector of parameters $\theta^0 \in\mathbb{R}^{I-1}$ such that
$\| \theta_h^0 \|$ is positive for all $h \in\mathcal{H}(\Delta_0)$
and zero otherwise. Having observed a contingency table $\mathbf{n}$, we
seek to recover $\Delta_0$. That is, our goal is to identify those
block components of $\theta^0$ having positive norms. To this end, we
define the group-lasso estimator for log-linear models to be the
solution of the concave optimization problem
%
%
%e10 ###
\begin{equation}\label{eqPl}
\max_{\theta\in\mathbb{R}^{I-1}} P_{\Lambda}(\theta)\equiv\max
_{\theta\in\mathbb{R}^{I-1}} \biggl\{\frac{1}{N} \ell(\theta) -
\lambda\sum_{h} \lambda_{h} \| \theta_h \|\biggr\} ,
\end{equation}
with $\ell(\cdot)$ defined as in (\ref{eqsat}) and $\Lambda= \{
\lambda, \{ \lambda_h, h \neq\varnothing\} \}$ a set of given
tuning parameters. The parameter $\lambda$ controls the overall effect
of the penalty and should be a function of the sample size, while the
block parameters $\lambda_h$ allows for specific penalties depending on
the sizes of the individual blocks.
A reasonable choice for these tuning parameters is $\lambda_h = \sqrt
{d_h}$, so that each block of coefficients is penalized proportionally
to its dimension, with larger blocks penalized more heavily.

%where $\ell(\cdot)$ is defined by (\ref{eqsat}) and $\m{pen}(\cdot)$
%assigns a penalty to every block $\theta_h$ that is non-zero. Ideally,
%the function $\m{pen}$ should satisfy two requirements. First, it
%should act as a thresholding function by either keeping or killing
%(i.e. setting to zero) each block $\theta_h$, $h \subset\mathcal{K}$,
%$h \neq\varnothing$. Secondly, it should be reasonably well behaved
%(e.g. convex) so that the problem (\ref{eqpeneq}) is computationally
%feasible.

The group-lasso penalty appearing in (\ref{eqPl}) was first proposed
by Yuan and Lin \cite{YL06} in the context of linear Gaussian models under ANOVA
settings (see also Nardi and Rinaldo~\cite{NR08}). It is specifically designed to
produce sparsity in the vector of estimated coefficients at the block
level. It is obtained as compositions of the $\ell_1$ norm over
quadratic norms of the individual blocks. The quadratic norms of
individual blocks promote non-sparsity, whereas the $\ell_1$ norm
applied to the resulting block norms, promotes block sparsity.
The group-lasso methodology of Yuan and Lin~\cite{YL06} was further extended to
logistic regression models by Meier, van der Geer and B\"{u}hlmann~\cite{MVB06} and to log-linear models by
Dahinden et al. \cite{LOGLASSO06}, which inspired our work.

\begin{lemma}\label{lamexistsol}
The vector $\widehat{\theta} \in\mathbb{R}^{I-1}$ is an optimizer of
\eqref{eqPl} if and only if there exists a~vector $\widehat{\eta}
\in
\mathbb{R}^{I-1}$ such that, for any $h$,
%
%
%e11 ###
\begin{equation}\label{eqoptim}
-\frac{1}{N}\mathrm{U}^{\top}_h (\mathbf{n} - \widehat{\mathbf
{m}}) + \lambda
\lambda_h \widehat{\eta}_h = 0,%\frac{\widehat{\theta}_h}{\Vert \widehat{
\end{equation}
where
\[
\widehat{\eta} = \cases{
\displaystyle\frac{\widehat{\theta}_h}{\Vert \widehat{\theta}_h\Vert _2}, & \quad$\mbox{if }
\widehat{\theta}_h \neq0$,\vspace*{2pt}\cr
\lambda_h \widehat{z}_h, & \quad$\mbox{if } \widehat{\theta}_h = 0$}
\]
with $\| \widehat{z}_h \| \leq1$, and $\widehat{\mathbf{m}} = \frac
{N}{\langle\exp(\vphantom{\tfrac{a}{a}}\mathrm{U} \widehat{\theta}),\mathbf{1}\rangle}
\exp(\mathrm{U}
\widehat{\theta})$. The solution is unique if $\| \widehat{z}_h \| < 1$
for each~$h$ for which $\widehat{\theta}_h = 0$.
%admits a unique optimizer whose $h$-block component satisfies
%-\frac{1}{N}\m{U}^{\top}_h (\mathbf{n} - \widehat{\mathbf{m}}) +
% \mbox{if } \widehat{\theta}_h \neq0\\
%-\frac{1}{N}\m{U}^{\top}_h (\mathbf{n} - \widehat{\mathbf{m}}) +
%where the vectors $\{ \widehat{z}_h \}$ satisfy $\| \widehat{z}_h \|
%If
% = \mathbb{E}_{\widehat{\theta}} \mathbf{n}$.
\end{lemma}

Having obtained the group-lasso estimator $\widehat{\theta}$, the model
selection step entails building an estimate of the true model $\Delta
_0$ by extracting the blocks of $\widehat{\theta}$ with positive norm
and then build a hierarchical model $\widehat{\Delta}$ as illustrated
in Table \ref{tabalgo}.
%One may say that this procedure is effective at recovering the
%underlying set of interactions if, with high probability, $\widehat{
%model selection consistency. Notice that, since we are only concerned
%with finding good estimators of $\Delta_0$, it is not required of $
%of $\theta^0$, besides the ones leading to model selection
%consistency. In fact, we will see in Section \ref{secclt} that $
There are two advantages in using the group-lasso estimator for
estimating $\Delta_0$ rather than traditional methods of model
selection based on sequential testing of a potentially very large
number of competing models.
The first advantage is that the methodology described in Table \ref
{tabalgo} only involves determining a penalized maximum likelihood
estimator of $\theta^0$ and thus requires solving only one convex
optimization problem (albeit a hard one, see the discussion below). In
contrast, classical model selection procedure requires fitting and
comparing a~number of different models which, even for tables with a
small number of variables, can be unfeasible. %On the other hand, the
%group-lasso estimator, being a penalized maximum likelihood estimator
%of the saturated model, has largest complexity over all competing
%sub-models, so that it would not be feasible to compute for tables so
%large that only smaller models can be fitted.
The second advantage is that the group-lasso estimator always return
a~model for which the maximum likelihood estimate exists, a fact that is
not guaranteed by the computational procedures currently used in
practice. See Fienberg and Rinaldo~\cite{FIER06} for more details.

%
%t1 ###
\begin{table}
\caption{The group-lasso model selection for hierarchical log-linear models}
\label{tabalgo}
\begin{tabular*}{210pt}{@{\extracolsep{\fill}}l@{}}
\hline
1. Obtain the log-linear group-lasso estimator,\\
\multicolumn{1}{c}{$
\widehat{\theta} = \operatorname{argmax}_{\theta\in\mathbb{R}^{I-1}}
P_{\Lambda}(\theta).
$}\\
2. Extract the set of non-zero blocks from $\widehat{\theta}$,\\
\multicolumn{1}{c}{$
\widehat{\mathcal{H}} = \{ h \dvt\| \widehat{\theta}_h \| > 0 \}.
$}\\
3. Recover the hierarchical log-linear model from $\widehat
{\mathcal{H}}$,\\
\multicolumn{1}{c}{$
\widehat{\Delta} = \{ h' \dvt h' \subseteq h, \mbox{ for some } h
\in
\widehat{\mathcal{H}} \}.
$}\\
\hline
\end{tabular*}
\end{table}

\begin{remark*}
Though motivated by model selection with hierarchical models, our
analysis below will actually show the log-linear group lasso estimator
can asymptotically recover any log-linear interaction subspace. This is
the reason why in equation \eqref{eqH} we used the more general
notation~$\mathcal{H}$ to encode the interactions of a hierarchial
log-linear model $\Delta$. Thus, the last step in the algorithm
described in Table \ref{tabalgo}, which forces the estimated model to
be hierarchical, should not be used if interested in general
interaction models. Furthermore, while asymptotically this step is
unnecessary for a hierarchical model, we nonetheless believe it would
improve the finite sample performance of the algorithm.
\end{remark*}

%One could make different choices for the group penalty function. For
%instance, \cite{TROPP05} considers a penalty term which is built as
%the composition of the $\ell_1$ norm over the $\ell_{\infty}$ norms of
%individual blocks. This particular choice would assign a milder
%penalty for complexity, for the same set of tuning parameters $

\subsection*{Model complexity and computational considerations}

We now make some remarks on the complexity of the class of log-linear
models. The model selection problem for log-linear models is
characterized by a combinatorial explosion in the number of possible
models that is much larger than for linear and many generalized-linear
models. Combined with the fact that the notion of sample size is also
quite different, as it refers to total number of counts $N$ in our
settings and to the total number of observed counts $\prod_k I_k$ in
the linear and generalized linear model settings, a~direct comparison
between the computational burden of our estimator versus the lasso or
group-lasso procedure for those models is not entirely adequate.

%a more complex problem than model selection for linear and many
%generalized-linear models.

It follows from our combinatorial definition of a general log-linear
model from Section~\ref{secloglinsubs} that, for a $K$-way table,
log-linear models can be represented as subsets of the set of all the
$2^K - 1$ possible interactions, except the one specified by the empty
set. Table \ref{tabnK} displays the number of log-linear models of
different types as a function of $K$. For a given $K$, the number of
possible unrestricted log-linear model is $2^{2^K-1}$, which
super-exponential in $K$. The number of hierarchical models, for which
no closed form expression is available, is much smaller, but still
appears to grow extremely fast (roughly exponentially) in $K$. For the
even smaller subclasses of graphical\footnote{The number of possible
graphical models is $\sum_{i=0}^K{K \choose i}2^{{i \choose2}}$.} and
decomposable models, the space of possible models is still extremely
large, for small values of $K$. For example, for a binary 5-dimensional
table with a total of $32$ observed counts, there are $1233$
decomposable models!
Due to to this tremendous combinatorial complexity, model selection by\vadjust{\goodbreak}
exhaustive model search is computationally prohibitive and, in fact,
unfeasible even for very small tables.
Indeed, model selection techniques for log-linear models must rely
necessarily on very greedy algorithms and are often designed to
consider only graphical or decomposable models. See, for instance, the
Baysian procedure of Dobra and Massam \cite{DM10}, and the frequentis model selection
search based on asymptotic $\chi^2$ testing implemented in the software
\texttt{MIM}, described in Edwards \cite{EDWARDS00}.

In contrast, the group-lasso procedure we consider is parametrized, in
both the likelihood and the penalty part, by only $2^K - 1$ terms,
which correspond to the direct sum decomposition of $\mathbb
{R}^{\mathcal{I}}$ into the orthogonal interaction subspaces. Although
this number grows exponentially in $K$, it is still order of magnitudes
smaller than the possible number of models. Indeed, based on Table \ref
{tabnK}, it appears to be loosely logarithmic in the number of
possible models. As a result, even though the computational complexity
of (\ref{eqPl}) may be high, it is significantly smaller than
exhaustive model selection and many greedy model selection algorithms.

%
%t2 ###
\begin{table}
\caption{Number of log-linear models of different types as a function
of $K$. Source: Lauritzen~\protect\cite{LAU02}}
\label{tabnK}
\begin{tabular*}{\textwidth}{@{\extracolsep{\fill}}llllll@{}}
\hline
& \multicolumn{5}{c@{}}{K}\\[-6pt]
& \multicolumn{5}{c@{}}{\hrulefill}\\
Type & 1 & 2 & 3 & 4 & 5\\
\hline
Unrestricted & 2 & 8 & 128 & 32,768 & 2,147,483,648\\
Hierarchical & 2 & 5 & \phantom{0}19 & 167 & 7580\\
Graphical & 2 & 5 & \phantom{0}18 & 113 & 1450\\
Decomposable & 2 & 5 & \phantom{0}18 & 110 & 1233\\
\hline
\end{tabular*}
\end{table}

%On the other hand the total number of cell is of order $O(I_{k_{

To our knowledge, two algorithms for computing the group-lasso
estimates are currently available: the block-coordinate descent method
of Meier, van der Geer and B\"{u}hlmann\cite{MVB06}, developed for the specific case of logistic
regression with grouped variables but immediately extendible to our
settings, and the path-following algorithm of Dahinden et al. \cite{LOGLASSO06} for
general log-linear models on binary variables. Both methods showed good
performance on simulated data and have been applied successfully to
real-life datasets, with the tuning parameters chosen by
cross-validation. Those results corroborate our theoretical findings
that group-lasso estimator possess good theoretical properties and is a
valuable alternative to greedy model selection procedures.
%Nonetheless we remark that a thorough study of the computational
%complexity of the group-lasso solution (\ref{eqPl}) is in order.
For completeness, we reference the more recent works by Roth and Fisher \cite{RF08},
Puig, Wiesel and Hero \cite{PUIG}, Yuan, Roshan, and Zou \cite{Y09} and Friedman, Hastie and Tibshirani \cite{FHT10}.

Nonetheless, we remark that the development of efficient computational
methods for calculating the group-lasso solution (for log-linear as
well as for linear model) with proven performance particularly in very
high-dimensional settings still remains an open problem.

\section{Asymptotic analysis}\label{secasymptotics}

%In the remainder of the paper we will perform a `large $p$ and large
%$N$' type of asymptotic analysis of the model selection procedure
%described in Table \ref{tabalgo} and of the properties of the
%group-lasso estimator.

In this section, we provide the main results of the paper. We perform
here a `large $p$ and large~$N$' type of an asymptotic\vadjust{\goodbreak} analysis of the
model selection procedure described in Table \ref{tabalgo} and of the
properties of the group-lasso estimator. We consider a rather general
double-asymptotic framework, in which we allow both the sample size and
the complexity of the statistical model to grow simultaneously.
In particular, we assume a sequence of statistical experiments
consisting of log-linear models over an increasingly large set of cell
combinations, implied by both a growing number of categorical variables
and a growing number of levels for the variables, and with increasing
sample size.
To formally represent this sequence of experiments, we will introduce a
`time' variable $n$, which serves merely as an index and is not
necessarily a quantification of the rate of increase of the sample
size. Intuitively, the larger the index $n$, the bigger the contingency
table, the larger the sample size and the more complex the model
selection problem.

To be specific, at time $n$,
\begin{itemize}
\item it is available a multinomial sample of size $N_n$ from the joint
distribution of $K_n$ categorical variables, each defined over a finite
set $\mathcal{I}_{k_n} = \{1, \ldots, I_{k_n} \}$, $k_n =1 ,\ldots,
K_n$; the support of this distribution is the set $\mathcal{I}_n =
\bigotimes_{k_n} \mathcal{I}_{k_n}$ of all cell combinations, of
cardinality $I_n = \prod_ {k_n} I_{k_n}$;
\item the true underlying distribution is defined by a hierarchical
log-linear model $\Delta_n$, as described in Section \ref
{secloglinsubs}: the observed cell counts come from an exponential
family distributions with log-likelihood function (\ref{eqsat}) and
true natural parameter $\theta^0_n \in\mathbb{R}^{I_n - 1}$, such that
$\| \theta^0_{h_n} \| > 0 $ for $h_n \in\mathcal{H}_n$ and $\|
\theta
^0_{h_n} \| = 0$ for $h_n \notin\mathcal{H}_n$, with $\mathcal{H}_n$
defined as in~(\ref{eqH}); the corresponding vector of cell
probabilities is denoted with $\bolds{\pi}_n^0 \in\mathbb
{R}^{\mathcal
{I}_n}$ and the mean vector $N_n \bolds{\pi}_n^0$ with $\mathbf{m}_n^0$;
\item the vector of true parameters $\theta^0_n$ is estimated by
solving the program (\ref{eqPl}) with tuning parameters $\Lambda_n =
\{ \lambda_n, \{ \lambda_{h_n}, h_n \neq\varnothing\} \}$;
\item the group-lasso estimate $\widehat{\theta}_n$ is then used to
estimate $\Delta_n$ as described in Table \ref{tabalgo}, leading to
the optimal selected model $\widehat{\Delta}_n$.
\end{itemize}
In the rest of the article, we will use the notation $\{t_n \} \in
\bigotimes\mathbb{R}^{k_n}$ to denote a sequence of vectors such that
$t_n \in\mathbb{R}^{k_n}$, for every $n$.

We remark that the true model at each `time point' $n$ needs not be
related with the true models at different values of $n$. The sequential
setting we adopt is a convenient device for representing very generally
an asymptotic framework for log-linear model selection with a diverging
number of parameters; in fact, there are many factors that may increase
the complexity of a log-linear model (e.g., number of variables, number
of interactions in the model, number of levels for each variable) that
we found it convenient to just allow each of them to change at every $n$.

In our sequential setting, the probability spaces are allowed to change
with $n$ and, when we speak of convergence in probability to a constant
or of tightness with respect to the index $n$, we explicitly refer to a
sequence of different probability measures. Accordingly, we will use
the stochastic small and large order notation ${o}_{P_n}$ and ${O}_{P_n}$
respectively with an index $n$ for the probability measures. This
notation is well defined: see, for instance, Schervish (\cite{MARK98}, Definition~7.11
and Lemma 7.12).

Projecting down the true parameter $\theta^0_n$ into $\mathbb
{R}^{d_{\mathcal{H}_n}}$ we write $\theta^0_{\mathcal{H}_n}$. %We
%will
%indicate with $\{ \mbf{\pi}_n^0 \}_n$ the sequence of true probability
%vectors and with $\{ \mathbf{m}_n^0 \}_n$ the sequence of mean vectors,
%with $\mathbf{m}_n^0 = N_n \mbf{\pi}_n^0$, for each $n$.
One may take note that, for a single observation, the Fisher
information matrix at $\theta^0_{\mathcal{H}_n}$ is
\[
\mathrm{F}_{\mathcal{H}_n} = \mathrm{U}^{\top}_{\mathcal{H}_n} \bigl(
\mathrm{D}_{ \bolds{\pi}_n^0} - \bolds{\pi}_n^0 ( \bolds{\pi
}_n^0)^{\top} \bigr) \mathrm{U}_{\mathcal{H}_n},
\]
with maximal and minimal eigenvalues denoted by $l^{\max}_n$ and
$l^{\min}_n$, respectively.
The negative Hessian of the log-likelihood function is
%
%
%e12 ###
\begin{equation}\label{eqS}
\Sigma_{\mathcal{H}_n} = \mathrm{U}_{\mathcal{H}_n}^{\top} \biggl(
\mathrm{D}_{\mathbf{m}^0_n} - \frac{\mathbf{m}^0_n (\mathbf
{m}^0_n)^{\top}}{N_n}
\biggr) \mathrm{U}_{\mathcal{H}_n} = N_n \mathrm{F}_{\mathcal{H}_n},
\end{equation}
which is the covariance matrix of the natural sufficient statistics
$\mathrm{U}_{\mathcal{H}_n}^{\top} \mathbf{n}_n$ and the Fisher
information based
on a i.i.d. sample of size $N$.

In the reminder of the article, we will study some of the asymptotic
properties of the sequence of group-lasso estimates $\{ \widehat
{\theta
}_n \}_n$ generated according to the previous scheme. %Although framed
%in a more general setting, our results carry over to standard
%asymptotic framework in which only the sample size $N_n$ and the set
%of penalty parameters $\Lambda_n$ change with $n$.
%In our analysis we will establish a series of progressively stronger
%results, which, naturally, demand increasingly stronger assumptions.
In Section \ref{secmodselconsistency}, we prove the that the
group-lasso estimator is model selection consistent, that is,
%
%
%e13 ###
\begin{equation}\label{eqlimresult}
\lim_n \mathbb{P} ( \widehat{\Delta}_n = \Delta_n ) = 1.
\end{equation}
Finally, in Section \ref{secclt} we give a central limit theorem for
$\widehat{\theta}_n$.

%Models}

%s4.1 ###
\subsection{Model selection consistency}\label{secmodselconsistency}

Here, we derive sufficient conditions for the property of model
selection consistency (\ref{eqlimresult}). Our method of analysis is
based on linearizing the sub-gradient optimality conditions~(\ref
{eqoptim}) via a Taylor expansion around the sequence of true
parameters $\theta^0_n$. As it turns out, norm (or $l_2$) consistency
is necessary to guarantee enough stochastic control over the remainder
term of that expansion. To that end, we first establish norm
consistency (in Lemma \ref{lemestimationconsistecy}) concerning a
related optimization problem (see (\ref{eqtildetheta}) below).
The conditions we develop for model selection consistency are quite
similar in spirit to the ones arising from the study of sparse recovery
of a linear signals under Gaussian or white noise using the lasso
penalty (see, in particular, Wainwright \cite{W06}, Zhao and Yu \cite{ZY06}).

Recall the definition of $\mathcal{H}_n$ from (\ref{eqH}) and let
$\mathcal{H}^c_n = 2^{\mathcal{K}_n} \setminus( \mathcal{H}_n
\cup\varnothing)$, so that $\| \theta^0_{h_n} \| > 0$ for each
$h_n \in\mathcal{H}_n$ and $\| \theta^0_{w_n} \| = 0$ for each $w_n
\in\mathcal{H}^c_n$.
Consider the sequence of events indexed by $n$
%
%
%e14 ###
\begin{equation}\label{eqOn}
\mathcal{O}_n = \{ \| \widehat{\theta}_{h_n} \| > 0, \forall h_n
\in\mathcal{H}_n \} \cap\{ \| \widehat{\theta}_{w_n} \| =
0, \forall w_n \in\mathcal{H}^c_n \}.
\end{equation}
Then, the model selection consistency property (\ref{eqlimresult}) of
the group-lasso solutions $ \widehat{\theta}_n$ is equivalent to
convergence in probability of $\mathcal{O}_n$, namely $\lim_n \mathbb
{P} ( \mathcal{O}_n ) = 1$.
This in turn occurs if and only if
%
%
%e15 ###
\begin{equation}\label{eqOn1}
\lim_n \mathbb{P} ( \| \widehat{\theta}_{h_n} \| > 0, \forall h_n
\in\mathcal{H}_n ) =1
\end{equation}
and
%
%
%e16 ###
\begin{equation}\label{eqOn2}
\lim_n \mathbb{P} ( \| \widehat{\theta}_{w_n} \| = 0, \forall w_n
\in\mathcal{H}^c_n ) =1.
\end{equation}

In this section, we will provide sufficient conditions for \eqref
{eqOn1} and \eqref{eqOn2}.

Our method of analysis relies on the primal-dual witness construction
of Wainwright~\cite{W06}, which we summarize below.
\begin{enumerate}[1.]
\item Solve a restricted group-lasso problem
%
%
%e17 ###
\begin{equation}\label{eqtildetheta}
\widetilde{\theta}_n = \mathop{\operatorname{argmax}}_{\theta\in\mathbb
{R}^{d_{\mathcal
{H}_n}}} P_{\Lambda}(\theta),
\end{equation}
where $\Lambda= \{ \lambda_n, \lambda_h, h \in\mathcal{H}_n(\Delta
_0)\}$.
\item Choose a vector $\widehat{\eta}_{\mathcal{H}_n}$ that belongs to
the subdifferential of group lasso penalty $\lambda_n \sum_{h, h \in
\mathcal{H}_n} \lambda_h \| x_h \|$ evaluated at $\widetilde{\theta
}_n$, and solve for a $(I-1 - d_{\mathcal{H}_n})$-dimensional vector
\[
\widehat{\eta}_{\mathcal{H}_n^c} = \operatorname{vec} \{ \lambda_{w_n}
\widehat
{z}_{w_n}, w_n \in\mathcal{H}_n^c \}
\]
the optimality conditions \eqref{eqoptim}. Check that $\|\widehat
{z}_{w_n}\|< 1$, for all $w_n \in\mathcal{H}_n^c $.
\item The vector $\hat{\theta}_n = \operatorname{vec}\{\widetilde{\theta
}_n,0\}
\in\mathbb{R}^{I-1}$ is the unique solution (see Lemma \ref
{lamexistsol}) of \eqref{eqoptim} and the event in equation \eqref
{eqOn} holds.
\end{enumerate}

Since the optimality conditions \eqref{eqoptim} are non-linear
functions of $\theta\in\mathbb{R}^{I-1}$, step 2. entails first a
linearization step to bound the difference between $\widetilde{\theta
_n}$ and $\theta^0_{\mathcal{H}_n}$ by a Taylor series expansion, and
then showing that the resulting remainder term vanished in probability.
%%is of order ${o}_{P_n^0}(\| \widetilde{\theta}_n - \theta^0_n\|
%)$.
This, in turn, follows from the next result, which established that
$\widetilde{\theta_n}$ is a norm consistent estimator of $\theta
^0_{\mathcal{H}_n}$.

\begin{lemma}\label{lemestimationconsistecy}
Assume
\begin{enumerate}[1.]
\item$\mathrm{[NC.1]}$
$d_{\mathcal{H}_n} ={o}(N_n)$;
\item$\mathrm{[NC.2]}$ $0 < D_{\min} < l^{\min}_n \leq l^{\max}_n
< D_{\max
} < \infty$, here $l^{\min}_n$, $l^{\max}_n$ are the minimal and
maximal eigenvalue of the Fisher information matrix, respectively;
%$\frac{\sqrt{l^{\max}_n}}{l^{\min}_n} = O(1)$
%where $l^{\min}_n$, $l^{\max}_n$ are the minimal and maximal
%eigenvalue of the Fisher information matrix, respectively;
%
\item$\mathrm{[NC.3]}$ for any $D>0$,
%
%
%e18 ###
\begin{equation}\label{eqNCeigen}
\sup\Biggl\{ | \mathbb{E}_{\theta_n}[ \langle a, \mathrm{U}_{\mathcal
{H}_n}^\top X \rangle] | \dvt\| \theta_{n}
- \theta^0_{\mathcal{H}_n} \| \leq D \sqrt{\frac{d_{\mathcal
{H}_n}}{N_n}}, \| a\| = 1 \Biggr\} = {O} \Biggl(\sqrt{\frac
{N_n}{d_{\mathcal{H}_n}}}\Biggr),
\end{equation}
where, for $\theta_n \in\mathbb{R}^{d_{\mathcal{H}_n}}$, $\mathbb
{E}_{\theta_n}$ denotes the expectation operator with respect to the
distribution $\operatorname{Multinomial}(1,\bolds{\pi}_n)$, with
\[
\bolds{\pi}_n = \frac{\exp(\mathrm{U}_{\mathcal{H}_n} \theta
_n)}{\langle\exp
(\mathrm{U}_{\mathcal{H}_n} \theta_n), \mathbf{1}\rangle};
\]
\item$\mathrm{[NC.4]}$
$\lambda_n = {O}( \frac{1}{\sum_{h_n \in\mathcal{H}_n} \lambda
_{h_n} } ) $.
%{O} ( \sqrt{\frac{d_{\mathcal{H}_n}}{N_n}} \frac{1}{\sum_{h_n \in
\end{enumerate}
Then,
%
%
%e19 ###
\begin{equation}\label{eqOpconsist}
\| \widetilde{\theta}_n - \theta^0_{\mathcal{H}_n} \| = {O}_{P_{n}^0}
\Biggl( \sqrt{\frac{ d_{\mathcal{H}_n} } {N_n}} \Biggr) = {o}_{P_n^0}(1).
\end{equation}
\end{lemma}

%The proof of Lemma \ref{lemestimationconsistecy} relies on a more
%general result about sequences of regular statistical experiments with
%diverging number of parameters, which is of independent interest and
%is described in Appendix B.

The proof of Lemma \ref{lemestimationconsistecy} relies on results of
Portnoy \cite{MR924876}. Indeed, conditions [NC.1] and~[NC.2] are essentially
derived from Portnoy \cite{MR924876}, Theorem 2.1. The more technical
condition [NC.3], also derived from Portnoy \cite{MR924876}, equation (2.4), is
needed to establish some control the order of magnitude of the
remainder term in the local quadratic approximation of the
log-likelihood function around $\theta^0_{\mathcal{H}_n}$, uniformly
over compact neighborhoods (see also Ghosal \cite{GHOSAL00}, for similar conditions).

Armed with \eqref{eqOpconsist}, we now proceed to prove the property
of model selection consistency for the group-lasso.

\begin{theorem}\label{thmmodselconsistency}
Assume the conditions of Lemma \ref{lemestimationconsistecy}.
Then, %of norm consistency [NC.1], [NC.3], [NC.3'] and [NC.4]. Then,
equation \eqref{eqOn1} holds if %(\ref{eqlimtheta}) holds if
%) $, assume
% \mbox{for all } n.
%
\begin{itemize}
\item$\mathrm{[MSC.1]}$ %A.3
letting $\alpha_n = \min_{h_n \in\mathcal{H}_n} \| \theta^0_{h_n}
\| $,
%
%
%e20 ###
\begin{equation}\label{eqconsist1}
\frac{1}{\alpha_n} \Biggl( \sqrt{\frac{d_{\mathcal{H}_n}}{N_n}} +
\lambda
_n \sqrt{\sum_{h_n \in\mathcal{H}_n} \lambda^2_{h_n}} \Biggr)
\rightarrow0,
\end{equation}
which, for $\lambda_{h_n} = \sqrt{d_{h_n}}$, simplifies to $\frac
{\sqrt
{d_{\mathcal{H}_n}}}{\alpha_n} ( \sqrt{\frac{1}{N_n}} + \lambda_n
) \rightarrow0$.
\end{itemize}
Equation \eqref{eqOn2} holds if %Equation (\ref{eqlimz}) holds if
\begin{itemize}
\item$\mathrm{[MSC.2]}$ (`almost' parameter orthogonality) for some
$\varepsilon\in(0,1)$ and for each $w_n \in\mathcal{H}_n^c$,
%[A.4]
%
%
%e21 ###
\begin{equation}\label{eqparort2}
\biggl\| \mathrm{U}^{\top}_{w_n} \biggl( D_{\mathbf{m}^0_n} -
\frac{\mathbf{m}^0_n (\mathbf{m}^0_n)^{\top}}{N_n}\biggr)
\mathrm{U}_{\mathcal{H}_n}
\Sigma_{\mathcal{H}_n}^{-1}
%( \m{U}^{\top}_{\mathcal{H}_n}
%( D_{\mathbf{m}^0_n} - \frac{\mathbf{m}^0_n(\mathbf{m}^0_n)^{
%)
\biggr\| < \frac{(1 - \varepsilon)}{|\mathcal{H}^c_n|};
\end{equation}
%
%(\mathbf{m}^0_n)^{\top}}{N_n}) \m{U}_{H_n} ( \m{U}^{
%
\item$\mathrm{[MSC.3]}$
%[A.5] (sparsity assumption)
%
\[
\lim_n \frac{|\mathcal{H}_n| \max_{h_n \in\mathcal{H}_n} \lambda
_{h_n}}{|\mathcal{H}_n^c| \min_{w_n \in\mathcal{H}_n^c}\lambda_{w_n}}
\leq1;
\]
\item$\mathrm{[MSC.4]}$%[A.6]
\[
\biggl( \max_{w_n \in\mathcal{H}_n^c} \frac{d_{w_n} }{\lambda_{w_n}^2}
\biggr)\frac{\log(I_n -1 - d_{\mathcal{H}_n})}{N \lambda_n^2 }
\rightarrow0, % \mbox{ and } \frac{
\]
which, for the choice $\lambda_{h_n} = \sqrt{d_{h_n}}$, becomes $
\frac
{\log(I_n -1 - d_{\mathcal{H}_n})}{N \lambda_n^2} \rightarrow0,
\rightarrow\infty$.

%) \rightarrow\infty,
%for each $w_n \in\mathcal{H}^c_n$,
%( \lambda_n \lambda_{w_n})^2 \frac{N_n}{d_{H_n}}
%( \lambda_n \lambda_{w_n})^2 \frac{N_n}{d_{H_n}}
%( \lambda_n \lambda_{w_n})^2 \frac{N_n}{d_{w_n}}
%which, for the choice $\lambda_{h_n} = \sqrt{d_{h_n}}$, requires $
\end{itemize}

\end{theorem}
%

%We point out that equation (\ref{eqlimz}) in particular demands
%asymptotic control over the some of the critical points of group-lasso
%program (\ref{eqPl}), which are required lie in the $\ell_2$--unit
%balls. In contrast, for the model selection consistency of the usual
%lasso procedure it is enough to consider only the $\ell_{
% Since in high dimensions the $\ell_{\infty}$ ball is much larger than
%the $\ell_2$ ball, one may be led to think that the sub-gradient
%optimality conditions for the group-lasso penalty are too harsh. %,
%and thus this penalty is not appropriate for large dimensional
%settings.
%Fortunately, the block penalties $\{ \lambda_h, h \in\mathcal{H}\}$
%can be adequately chosen to contrast this effect and, in fact, to
%allow for the model selection consistency of the group-lasso solutions
%under increasing complexity of the model, similarly to the lasso case.
%See, in particular, condition [MSC.4] below.

Condition [MSC.2] implies that
\[
\|\mathrm{W}_n\|=
\biggl\| \mathrm{U}^{\top}_{\mathcal{H}^c_n} \biggl( D_{\mathbf{m}^0_n} -
\frac{\mathbf{m}^0_n (\mathbf{m}^0_n)^{\top}}{N_n}\biggr) \mathrm
{U}_{\mathcal{H}_n}
\Sigma_{\mathcal{H}_n}^{-1}
%( \m{U}^{\top}_{\mathcal{H}_n} ( D_{\mathbf{m}^0_n} -
%)^{-1}
\biggr\| < (1 - \varepsilon),
\]
which is the equivalent of the so-called \textit{irreducibility condition}
appearing in the growing lasso literature (e.g., Wainwright \cite{W06} and Zhao and Yu \cite
{ZY06}). %In fact, under the Gaussian homoschedastic ensemble
%settings, the irreducibility condition on the design matrices is
%precisely an `almost' parameter orthogonality condition on the Fisher
%information matrix for the saturated model. The notable difference in
%our [MSC.2] is the need to account for the number of blocks of
%parameters which are zero.
Condition [MSC.3] is a~sparsity condition and condition [MSC.4]
provides some information on the rates of increase for the dimensions
for the subspaces of interactions not included in the true models (see
equation (15) (a) in Wainwright~\cite{W06}). Inspection of the proof will reveal
that, using a~simpler argument based on Chebyshev's inequality only,
and assuming $\lambda_{h_n} = \sqrt{d_{h_n}}$ for all $h$, [MSC.4]
reduces to $\lambda_n^2 N_n \rightarrow\infty$, so that model
selection consistency is obtained under conditions that do no depend on $I_n$.

\subsection{A central limit theorem for the log-linear group-lasso
estimator}\label{secclt}
Our final results concern the large sample properties of the
distribution of lasso group estimates~$\{ \widehat{\theta}_n \}_n$. In
addition to the conditions guaranteeing both norm and model selection
consistency, we need to impose further restrictions guaranteeing some
form of asymptotic normality under our double asymptotic framework.
%In our final result, we derive the limiting distribution of the
%sequence of lasso estimates
%The general idea of the proof is restrict to the sequence of events $
The main rationale behind retaining the set of assumptions for
consistency is that they allow us to work only with the simpler and
well-behaved sequence of events $\mathcal{O}_n$ defined in (\ref
{eqOn}), which converges in probability. %For example, within $
%estimates has dimension $d_{\mathcal{H}_n}$.

%In fact, under $\mathcal{O}_n$, there is a rather convenient linear
%approximation, shown in Equation (\ref{eqtheta}), between $\widehat{
%estimates, and $\theta^0_{\mathcal{H}_n}$, the positive components of
%the true underlying parameter, with full asymptotic control over the
%remainder term. Consequently, within $\mathcal{O}_n$, a central limit
%type of results for $\widehat{\theta}_{\mathcal{H}_n}$ are readily
%available.
In general, asymptotic normality under the double asymptotic settings
obtains under stricter assumptions than under standard (i.e., with
fixed-dimensional parameter space) asymptotic problems. Below, we
provide a series of conditions, each providing a sense that for large
enough $n$, the group-lasso estimates (appropriately rescaled and
translated) are close to a standard Normal distribution. %We point out
%that only the first two results in the theorem statement amount to a
%full central limit statement, while the third result only offers
%necessary condition for asymptotic normality.

%In particular, we consider two alternative set of assumptions
%guaranteeing asymptotic normality, which correspond each to a
%different type of double asymptotic settings, depending on how small
%the dimension $d_{\mathcal{H}_n}$ of the underlying model needs to be
%with respect to the sample size $N_n$. In the more restrictive case,
%asymptotic normality follow immediately from the classic
%Lindberg--Feller conditions. In the other case, a more involved method
%of proof is required, which is based on a multivariate extension of
%some results by \cite{MR0370871}.
To state our result, we need to formulate some notation.
Let $\mathrm{J}^0_{\mathcal{H}_n}$ be a $d_{\mathcal{H}_n} \times
d_{\mathcal
{H}_n}$ block-diagonal matrix whose $h_n$-block is the $d_{h_n} \times
d_{h_n}$ matrix
\[
\lambda_{h_n} \frac{1}{\| \theta^0_{h_n}\|} \biggl( \mathrm{I}_{d_{h_n}} -
\frac{\theta^0_{h_n}}{\Vert \theta^0_{h_n}\Vert _2} \biggl( \frac{\theta
^0_{h_n}}{\Vert \theta^0_{h_n}\Vert _2} \biggr)^{\top} \biggr),
\]
with $h_n \in\mathcal{H}_n$ and with $\mathrm{I}_{d_{h_n}}$ denoting the
$d_{h_n}$-dimensional identity matrix. Below, $\mathrm{G}_n$ will
denote a
$k \times d_{\mathcal{H}_n}$ matrix,
where $k$ is an arbitrary fixed number, such that
%
%
%e22 ###
\begin{equation}\label{eqGn}
\lim_n \mathrm{G}_n \mathrm{G}_n^{\top} = \mathrm{G}
\end{equation}
for some $k \times k$ nonnegative and symmetric matrix $\mathrm{G}$.

\begin{theorem} \label{thmclt}
Assume the conditions for norm and model selection consistency and let
%
%
%e23 ###
\begin{equation}\label{eqclt}
%( \Sigma^{1/2}_{H_n} + \Sigma^{-1/2}_{H_n} N_n \lambda_n
%J^0_{H_n} ) ( \widehat{\theta}_n - \theta^0_n ) +
X_n = \sqrt{N_n} \mathrm{F}_{\mathcal{H}_n}^{-1/2} \bigl( ( \mathrm
{F}_{\mathcal{H}_n} + \lambda_n \mathrm{J}^0_{\mathcal{H}_n} ) (
\widehat{\theta}_n%{\mathcal{H}_n}
- \theta^0_{\mathcal{H}_n} ) + \lambda_n \eta^0_{\mathcal{H}_n} \bigr),
\end{equation}
where
\[
\eta^0_{\mathcal{H}_n} = \operatorname{vec} \biggl\{
\lambda_{h_n} \frac{\theta^0_{h_n}}{ \| \theta^0_{h_n}
\|}, h_n \in\mathcal{H}_n\biggr \}.
\]
\begin{enumerate}[1.]
\item For each sequence $\{ \mathrm{G}_n \}$ of $k \times d_{\mathcal
{H}_n}$ matrices satisfying (\ref{eqGn}), $\mathrm{G}_n X_n$ converges
weakly to the $N_{k}(0,\mathrm{G})$ distribution if either the
$\mathrm{[CLT.LF]}$ condition
\[
d_{\mathcal{H}_n} = {o}(N^{1/2}_n)
\]
or both the $\mathrm{[CLT.Ma]}$ condition
\[
d_{\mathcal{H}_n} = {o} ( N_n )
\]
and $\mathrm{[CLT.Mb]}$ condition
\[
\max_{i \in\mathcal{I}_n} \bolds{\pi}_i^0 = {O} \biggl( \frac{1}{\sqrt{N_n
d_{\mathcal{H}_n}}} \biggr).
\]
hold.
\item If the $\mathrm{[CLT.BE]}$ condition
\[
d_{\mathcal{H}_n} = {o}(N^{2/7}_n),
\]
holds, then
\[
\sup_{A_n} | \mathbb{P}(X_n \in A_n) - \mathbb{P}(Z_n \in A_n) |
\rightarrow0,
\]
where $Z_n$ has a $N_{d_{\mathcal{H}_n}}(0,\mathrm{I}_{d_{\mathcal{H}_n}})$
distribution, the supremum is
taken over the convex sets $A_n$ in $\mathbb{R}^{d_{\mathcal{H}_n}}$
and convergence occurs at the rate ${O}( \frac{d_{\mathcal
{H}_n}^{7/2}}{N_n})$.\vspace*{-3pt}
%$X_n$ and $Z_n$, respectively. Then, if the $\m{[CLT.LF]}$ condition
%%%or [CLT.Ma] and [CLT.Mb] hold.The [CLT.LF] condition is the [CLT.LF]
%condition holds
%d_{\mathcal{H}_n} = o(N^{1/3}_n)
%holds, for each $\varepsilon> 0$ and $T>0$ there exists a $n^0(

%)$;
%$D_1$.
%$\max_{i \in\mathcal{I}_n} \pi_i^0\leq D_1 I_n^{-1}$, for some
%positive constant $D_1$;
%$d_{\mathcal{H}_n}=o( I_n )$;
%$\lambda_n = {O} ( \frac{1}{\sqrt{d_{\mathcal{H}_n}}} )$;
%$d_{\mathcal{H}_n} = o ( \sqrt{N_n} )$.
%Furthermore, if
%where $\mathcal{T} = \{ \{ t_n \} \dvt\| t_n \| \rightarrow
%for each $\{ t_n \} \otimes\mathbb{R}^{k_n}$ and each $\varepsilon> 0$
%there exists a $N$, which depends only on $\varepsilon$ and the
%sequence $
%| \phi_n(t_n) - \psi_n(t_n) | < \varepsilon.
\end{enumerate}
\end{theorem}

%Let $\m{F}_{\widehat{\mathcal{H}}_n}$ and $\m{J}_{\widehat{
%are consistent estimators.

The theorem indicates that the group-lasso estimate is asymptotically
unbiased and inefficient.
In fact, equation (\ref{eqclt}) demonstrates that the asymptotic
behavior of the group-lasso estimator is affected by two terms. One is
the bias term $ \lambda_n \eta^0_{\mathcal{H}_n}$ which depends on the
gradient of the penalty function at the true parameter. The other term
$\mathrm{J}^0_{\mathcal{H}_n}$ is the Hessian at the true parameter
of the
penalty function, a positive definite matrix which inflates the inverse
Fisher information. Both these terms are asymptotically significant and
indicate that the group-lasso estimates may lack asymptotic optimality.
Note that this phenomenon is probably quite general (see also Fan and Peng \cite{fan-2004-32}, Theorem 2).

Both Condition [CLT.LF] and conditions [CLT.Ma] and [CLT.Mb] guarantee
the asymptotic normality of a \textit{fixed} number of linear combinations
of the coordinates of~$\widehat{\theta}_n$. In particular, it includes
that case of $
\mathrm{G}_n = [{\mathrm{I}_{k}\enskip \mathrm{O}}]$,
where $\mathrm{O}$ is a $k \times(d_{\mathcal{H}_n} - k)$ matrix of zeros.
For this choice, the marginal asymptotic normality of any fixed number
of coordinates of~$\widehat{\theta}_n$ is guaranteed. Condition
[CLT.LF] results from a simple Lindberg--Feller argument, whereas
conditions [CLT.Ma] and [CLT.Mb] follow by adapting and generalizing
some proofs in Morris \cite{MR0370871}. We note that [CLT.Mb] may be replaced
by $\max_{i \in\mathcal{I}_n} \bolds{\pi}_i^0\leq C I_n^{-1}$, for some
positive constant $C$ (see, e.g., Quine and Robinson \cite{QuinRobinorm1984}, Theorem
1). Then, in order for the theorem to hold, one
has to further assume that
$
I_n = {o}(\sqrt{d_{\mathcal{H}_n} N_n})
$,
which is compatible with the conditions for norm consistency.

Condition [CLT.BE] is a full central limit type of results for the
group-lasso estimator and is based on a multivariate Barry--Esseen type
of bound found in Bentkus \cite{BEN03}. As it is usual with uniform results of
this type, it is necessary to control the fluctuations of third order
moments, and, consequently, to have a rather large sample size. To our
knowledge, this is the best rate available. See also Portnoy \cite{PORTNOY86}
for a similar result requiring only a~rate $d^2_{\mathcal{H}_n} =
{o}(N_n^{1/2})$, whose applicability and relevance to our problem is
however unclear.\vspace*{-3pt}

\section{Conclusions}\label{secconclusions}\vspace*{-3pt}

In this article, we studied some asymptotic properties of the
group-lasso estimator. Our results show that this estimator can be used
to recover asymptotically the true underlying model under conditions
that allow for a model complexity increasing with the sample size and
also for a number of cells larger than $N$.

Our setting, analysis and results differ from existing analyses of
$\ell
_1$ regularized least square problems in a few aspects.
Firstly, unlike the case of regularized least squares or Gaussian error
problems, the first order optimality conditions for the group-lasso
program are nonlinear in the parameters. As model selection consistency
hinges upon establishing appropriate bounds for the norms of the
differences between the blocks of true and estimated parameters, our
strategy was to linearize the sub-gradient equations via a first order
Taylor expansion. This expansion, in turn, is valid provided one has
enough control over the remainder term, which we achieved by proving
the norm consistency property for the group-lasso estimate. Thus, in
our settings, norm consistency is necessary for model selection
consistency. In contrast, for quadratic problems, whose first order
conditions are linear in the parameters, norm consistency does not
appear to be needed, although, it may still be important for central
limit results, like in our case.

Secondly, we did not concern ourselves with any form of consistency
other than the model selection consistency. However, other forms of
consistency are also relevant for the class of models presented here.
In particular, we mention the general risk consistency and the $\ell_2$
consistency of the penalized estimators for generalized linear model
and logistic regression models by van de Geer \cite{sara106,sara206} and Meier,
van der Geer and B\"{u}hlmann~\cite{MVB06}, respectively, where non-asymptotic bounds and oracle
inequalities are available. Finally, a~rather general framework for
proving norm consistency of penalized maximum likelihood estimators
under decomposable regularizers which may be applicable to our problem
is presented in Negahban et al. \cite{NRWY10}.

Finally, in our problem we do not need to worry about random design.
%%This is a consequence of the contingency table settings and
%simplifies the analysis.
In fact, as we are working with exponential families of distribution,
the Fisher information matrix is data-independent. Consequently, unlike
for example the case of Gaussian ensembles, for model selection
consistency it is sufficient to impose analytic, and not stochastic,
conditions on the asymptotic behavior of the Fisher information. These
conditions (namely, the almost parameter orthogonality' condition
[MSC.2]) correspond to the various irreducibility condition used in the
lasso literature, that we equivalently formulate in terms of the Fisher
information.
%Except for the irreducibility conditions, we did not take advantage of
%the exponential nature of multinomial distribution, so that the norm
%consistency and central limit results of sections
%rather general properties that may hold for other families of
%distributions.

%Extension to exponential families.
%linearize the problem. thus we need control over the reminder term,
%which is guaranteed by the norm consistency.
%by the conditions on the Fisher info matrix

%s6 ###
\section{Proofs}\label{secproofs}

\begin{pf*}{Proof of Lemma \protect\ref{lemlik}}
We only provide a sketch of the proof and refer to Haberman~\cite{HAB74}, page~11, and Rinaldo \cite{ALEb06}, Lemma 2.2, for more details.
It is possible to show that the linear subspace $\mathcal{M} \cap
\mathcal{R}(\mathbf{1})^{\bot}$ and $\widetilde{\mathcal{M}}$ are
homeomorphic sets, and the one-to-one mapping between $\widetilde
{\bolds{\mu}} \in\widetilde{\mathcal{M}}$ and $\bolds{\beta} \in
\mathcal{M}
\cap\mathcal{R}(\mathbf{1})^{\bot}$ is given by
\[
\widetilde{\bolds{\mu}} = \bolds{\beta} + \mathbf{1} \log\biggl(\frac
{N}{\langle\exp(\bolds{\beta}) , \mathbf{1} \rangle} \biggr).
\]
Furthermore, for any $\mathbf{n}$, some algebra shows that
\[
\langle\mathbf{n}, \bolds{\beta} \rangle- N \log\langle\exp
(\bolds{\beta }), \mathbf{1} \rangle+ \log N! - \sum_{i \in
\mathcal{I}} \log\mathbf
{n}_{i}! = \ell^*(\widetilde{\bolds{\mu}}),\vspace*{-3pt}
\]
for each pair of homeomorphic points $\bolds{\beta}$ and
$\widetilde
{\bolds{\mu}}$. Then, for any full-rank matrix $\mathrm{U}$ with
$\mathcal{M}
\cap\mathcal{R}(\mathbf{1})^{\bot}$ as its column span, \eqref{eqexpfam}
and \eqref{eqhomeo} both follow from the previous displays by noting
that there exists also a one-to-one correspondence between $\mathcal{M}
\cap\mathcal{R}(\mathbf{1})^{\bot}$ and $\mathbb{R}^{I-1}$, given
by $\bolds{\beta} = \mathrm{U} \theta$.\vspace*{-4pt}
%Fortunately, through an appropriate re-parametrization, one can
%replace $\widetilde{\mathcal{M}}$ with the $(\mathrm{dim}(
%subspace of $\mathbb{R}^I$ spanned by $\mathbf{1}$. Specifically, for
%each $\mbf{\beta} \in\mathcal{M} \cap\mathcal{R}(\mathbf{1})^{\bot}$,
%set
%Then, %there exists a bijection between $\widetilde{\mathcal{M}}$ and $
%for each $\widetilde{\mbf{\mu}} \in\widetilde{\mathcal{M}}$ there
%exists one $\mbf{\beta} \in\mathcal{M} \cap\mathcal{R}(\mathbf{1})^{
%and, conversely, for each $\mbf{\beta} \in\mathcal{M} \cap
%][]{ALEb:06}
%Therefore, if $\m{U}$ is any full-rank matrix whose columns span $
%as
%%\begin{equation}\label{eqexpfam}
%%\end{equation}
\end{pf*}

\begin{pf*}{Proof of Lemma \protect\ref{lamexistsol}}
The first order optimality conditions for a vector
$\theta\in\mathbb{R}^{I-1}$ is $\mathbf{0} \in\partial P_{\Lambda
}(\theta) $, the
subdifferential set of $P_\Lambda(\theta)$.
The gradient of $\ell$ at a point $\theta\in\mathbb{R}^{I-1}$ is
%
%
%e24 ###
\begin{equation}\label{eqgrad}
\nabla\ell(\theta) = \mathrm{U}^{\top}\biggl ( \mathbf{n} - \biggl( \frac
{N}{\langle\mathbf{b},\mathbf{1}\rangle} \biggr) \mathbf{b} \biggr) = \mathrm
{U}^{\top} ( \mathbf{n} - \mathbf{m} ),\vspace*{-3pt}
\end{equation}
where $\mathbf{b} =\exp(\mathrm{U} \theta)$.
As for the penalty term, which is not differentiable when some of the
blocks are zero, standard subgradient calculus (see, e.g., Bertsekas \cite{BER95}) yields that for any $\theta\in\mathbb{R}^{I-1}$, the
subdifferential of the function $ x \mapsto\sum_{h} \lambda_{h} \| x_h
\|$ at $\theta$ is a subset of~$\mathbb{R}^{I-1}$ comprised by vectors
whose $h$-block component is
%
%
%e25 ###
\begin{equation}\label{eqsubpen}
\cases{
\lambda_h \displaystyle\frac{\theta_h}{\Vert \theta_h\Vert _2}, &\quad $\mbox{if } \theta_h
\neq0,$\vspace*{1pt}\cr
\lambda_h z_h, & \quad $\mbox{if } \theta_h = 0,$}\vspace*{-3pt}
\end{equation}
where $\| z_h \| \leq1$ for each $h$ such that $\theta_h = 0$.
%we use the fact that if $\Vert  \cdot\Vert $ is a norm in an Banach space
%$X$, then the subdifferential at a point $x$ is
% \mbox{if } x \neq0\\
%where $X^*$ denotes the dual space of $X$.
%Then, since the the space $L_2$ is self-dual, we conclude that for any
%$\theta\in\mathbb{R}^{I-1}$, the subdifferential of the function $ x
%is a subset of $\mathbb{R}^{I-1}$ comprised by vectors whose $h$ block
%component is
% \lambda_h z_h & \mbox{if } \theta_h = 0,\\
%.
%where $\| z_h \| \leq1$ for each $h$ such that $\theta_h = 0$.
Equations (\ref{eqgrad}) and (\ref{eqsubpen}) imply (\ref{eqoptim}).

As for uniqueness, we follow the proof of Lemma 2 in Wainwright \cite{W06}.
Suppose $\widehat{\theta}$ is an optimal solution to \eqref{eqPl}.
Then, by duality theory, given a subgradient $\widehat{\eta} \in
\mathbb
{R}^{I-1}$ any optimal solution $\breve{\theta}$ must satisfy
$\widehat
{\eta}^\top\breve{\theta} = \sum_{h} \lambda_h \| \breve{\theta
}_h \|$. This holds only if $\breve{\theta}_h = 0$ for all $h$ for which $\|
\widehat{\eta}_{h} \| < \lambda_h$. Thus, if there exists a solution
$\breve{\theta}$ to the problem \eqref{eqPl} different than
$\widehat
{\theta}$, it must satisfy $\breve{\theta}_h = 0$ for all $h$ such that
$\widehat{\theta}_h = 0$. Finally, uniqueness follows since $\ell$ is
strictly concave, though not necessarily strongly concave.\vspace*{-4pt} %By standard
%results, described for example in \cite{ERGM}, the solution to the
%likelihood equations is always unique, unless the maximum likelihood
%estimator does not exist, in which case for every sequence of
%parameters $\{ \theta_n\} _n$ such that
%$\| \theta_n \| \rightarrow\infty$. However, the penalty term would
%prevent this from happening. This, combined with the strict convexity
%of the $\ell_2$ norm, guarantees uniqueness.
\end{pf*}

\begin{pf*}{Proof of Lemma \protect\ref{lemestimationconsistecy}}
With some abuse of notation we write ${O}_{P^0_n}$ and ${o}_{P^0_n}$ to
refer to probabilistic statements for the sequence of probability
distributions indexed by $\{ \theta^0_{\mathcal{H}_n}\}_n$, with
$\theta
^0_{\mathcal{H}_n} \in\mathbb{R}^{d_{\mathcal{H}_n}}$ for each $n$.
We will first analyze the asymptotic behavior of $\ell(\theta
^n_{\mathcal{H}_n}) - \ell_{n}(\theta^0_{\mathcal{H}_n})$, uniformly
over sequences of the form $\theta^n_{\mathcal{H}_n} = \theta
^0_{\mathcal{H}_n} + \sqrt{\frac{d_{\mathcal{H}_n}}{N_n}} x_n $ with
$\| x_n \| \leq D $, for all $n$ and some $D>0$.
To this end, we follow the arguments used in Portnoy \cite{MR924876}, Theorem 2.4.
%Also, for simplicity, we multiply $P_{\Lambda_n}$ by $N_n$ in (
%satisfied for the sequence of exponential families under
%considerations, with $\eta_{\theta^0_n} = \nabla p_{\theta^0_n}(p_{
%= d_{\mathcal{H}_n}$. Then, for any finite $C>0$ and for each sequence
%$\{ x_n \}_n \in\bigotimes_{n} \mathbb{R}^{d_{\mathcal{H}_n}}$ with $
First off, notice that we can write
\[
\mathbf{n} = \sum_{j=1}^{N_n}{X_j},\vspace*{-3pt}
\]
where $X_1,\ldots,X_{N_n}$ are i.i.d. vectors in $\mathbb
{R}^{\mathcal{I}}$ distributed like a $\operatorname{Multinomial}(1,\bolds
{\pi }^0_n)$, with
\[
\bolds{\pi}^0_n = \frac{\exp(\mathrm{U}_{\mathcal{H}_n} \theta
^0_{\mathcal
{H}_n})}{\langle\exp(\mathrm{U}_{\mathcal{H}_n} \theta^0_{\mathcal{H}_n}),
\mathbf{1}\rangle}.\vadjust{\goodbreak}
\]
By a Taylor series expansion, the term $\ell_{n}( \theta
^0_{\mathcal{H}_n} + \sqrt{\frac{d_{\mathcal{H}_n}}{N_n}} x_n ) -
\ell_{n}(\theta^0_{\mathcal{H}_n})$ is equal to
%
%
%e26 ###
\begin{equation}\label{eqtaylor}
\hspace*{-5pt}\sqrt{\frac{d_{\mathcal{H}_n}}{N_n}} x_n^{\top} \mathrm
{U}_{\mathcal
{H}_n}^{\top} ( \mathbf{n}_n - \mathbf{m}^0_n ) - \frac{1}{2} \frac
{d_{\mathcal{H}_n}}{N_n} x_n^{\top} \Sigma_{\mathcal{H}_n} x_n +
\frac
{N_n}{6} \biggl( \frac{d_{\mathcal{H}_n}}{N_n} \biggr)^{3/2}\mathbb
{E}_{\theta^*_n}[ ( \langle x_n, \mathrm{U}_{\mathcal{H}_n}^\top
X_1 \rangle)^3], %\m{U}_n^{\top} ( \m{D}_{\mathbf{%m}^0_n} - \frac{
\end{equation}
where $\theta^*_n$ is on the line joining $\theta^n_{\mathcal{H}_n}$
and $\theta^0_{\mathcal{H}_n}$.
For the first term in \eqref{eqtaylor}, we have
\begin{eqnarray*}
\mathbb{E}\Vert  \mathrm{U}_{\mathcal{H}_n}^{\top} ( \mathbf{n}_n -
\mathbf
{m}^0_n ) \Vert ^2_2 & = & %\sum_{j=1}^{I_n} ( \mathbb{E}
%u_{l,j}u_{k,j} \mathbb{E} (n_l - m_l) (n_k - m_k) )\\
%& = & \sum_{i,j} \mathbf{u}_i^{\top} ( \m{D}_{\mathbf{m}^0_n} - \frac{{
%& = &
\operatorname{tr}\biggl( \mathrm{U}_{\mathcal{H}_n}^{\top} \biggl( \mathrm
{D}_{\mathbf{m}^0_n} -
\frac{\mathbf{m}^0_n (\mathbf{m}^0_n)^{\top}}{N_n} \biggr) \mathrm
{U}_{\mathcal
{H}_n} \biggr)\\
& \leq& N_n l_n^{\operatorname{max}} d_{\mathcal{H}_n}.\\
\end{eqnarray*}
%
%
%where $\lambda^n_{\m{max}} = l^n_{\m{max}} N_n$ denotes the largest
%eigenvalue of the Fisher information matrix.
Thus, by Markov and Cauchy--Schwarz inequalities,
%
%
%e27 ###
\begin{equation}\label{eqE1}
\Biggl| \sqrt{\frac{d_{\mathcal{H}_n}}{N_n}} x_n^{\top} \mathrm
{U}_{\mathcal
{H}_n}^{\top} ( \mathbf{n}_n - \mathbf{m}^0_n ) \Biggr| = {O}_{P_n^0}
\bigl( d_{\mathcal{H}_n} \sqrt{ l_n^{\operatorname{max}}} \| x_n \| \bigr) =
{O}_{P_n^0} ( d_{\mathcal{H}_n} \| x_n \| ), %\| x_n \| k_n
%{O}_{P_{\theta^0_{\mathcal{H}_n}}} ( \sqrt{l_n^{\m{max}}} ).
\end{equation}
where the last identity follows from the eigenvalue assumption [NC.2].
% \sqrt{\frac{k_n}{N_n}} x_n^{\top} \m{U}_n^{\top} ( \mathbf{n}_n - {
%& & &\\
%E_1 & + & E_2. &\\
Next, using the fact that the Fisher information matrix is positive
definite for each $n$, and once again~[NC.2],
%
%
%e28 ###
\begin{equation}\label{eqE2}
- \frac{1}{2} \frac{d_{\mathcal{H}_n}}{N_n} x_n^{\top} \Sigma
_{\mathcal
{H}_n} x_n \leq- \frac{1}{2} \frac{d_{\mathcal{H}_n}}{N_n} N_n \| x_n
\|^2 l^n_{\min} x \leq- \frac{1}{2} {O} ( d_{\mathcal{H}_n} \| x_n
\|^2 ),
\end{equation}
which bounds the second term \eqref{eqtaylor}. Finally, the assumption
[NC.3] yields (see Portnoy~\cite{MR924876}, Theorem 2.4)
\[
\frac{N_n}{6} \biggl( \frac{d_{\mathcal{H}_n}}{N_n} \biggr)^{3/2}\mathbb
{E}_{\theta^*_n}\biggl[ \biggl| \biggl\langle\frac{x_n}{\|x_n\|}, \mathrm
{U}_{\mathcal{H}_n}^\top X_1 \biggr\rangle\biggr|^3\biggr] = {O} (
d_{\mathcal{H}_n} \| x_n \| ).
\]

Combining the previous display with (\ref{eqE1}), (\ref{eqE2}) and
with \eqref{eqtaylor}, we obtain, by choosing~$D$ large enough that,
for each $\varepsilon> 0$, and all $n$ large enough
\[
\mathbb{P} \Biggl( \sup_{ \{ x_n, \| x_n \| = C \}} \ell
_{n}\Biggl( \theta^0_{n} + \sqrt{\frac{d_{\mathcal{H}_n}}{N_n}} x_n
\Biggr) - \ell_{n}(\theta^0_n) < 0 \Biggr) > 1- \varepsilon.
\]
The strict concavity of $\ell$ (warranted by [NC.2]) further guarantees
that, for all $n$ large enough, with probability tending to one, there
are no other maximizers of $\ell$ outside the ball $\{ \theta
^n_{\mathcal{H}_n} \dvt\theta^n_{\mathcal{H}_n} = \theta
^0_{\mathcal
{H}_n} + \sqrt{\frac{d_{\mathcal{H}_n}}{N_n}} x_n \}$.

Next, we consider the difference in the penalty terms between $\theta
^0_{\mathcal{H}_n}$ and $\theta^n_{\mathcal{H}_n} \equiv\theta
^0_{\mathcal{H}_n} + \sqrt{\frac{d_{\mathcal{H}_n}}{N_n}} x_n$. %,
%with
%$\| x_n \|\leq C$.
By a first order Taylor expansion,
\[
N_n \lambda_n \biggl( \sum_{h \in\mathcal{H}_n} \lambda_{h} \Vert  \theta
^n_h\Vert _2 - \sum_{h \in\mathcal{H}_n} \lambda_{h} \Vert \theta^0_{h}\Vert _2
\biggr)  = %\geq& N_n \lambda_n \sum_{h_n \mathcal{H}_n}
%- \| \theta^0_{h_n} \| )\\
%& = &
 N_n \lambda_n \sum_{h_n \in\mathcal{H}_n} \lambda_{h_n} \sqrt
{\frac
{d_{\mathcal{H}_n}}{N_n}} (x^*_{h_n})^{\top} \frac{\theta
^0_{h_n}}{\|
\theta^0_{h_n}\|},
\]
where $x^*_n$ lies between $0$ and $x_n$. Using the Cauchy--Schwarz
inequality, the absolute value of the last quantity is bounded by
\[
\lambda_n \sqrt{N_n d_{\mathcal{H}_n}} \| x_n \| \biggl( \sum_{h \in
\mathcal{H}_n} \lambda_h \biggr).\vspace*{-2pt}
\]

%By the triangle inequality,
%N_n \lambda_n ( \sum_{h \subseteq\mathcal{K}_n} \lambda_{h} \Vert
%where the summations are taken over set $h \neq\varnothing$. Then, in
%order for the penalty component to be of smaller order of magnitude
%than the other terms, we must have
%%\lambda_n = {O} ( \sqrt{\frac{I_n - 1}{N_n}} \frac{1}{\sum_{h

Under the assumed conditions, we see that $\| \widetilde{\theta}_n -
\theta^0_{\mathcal{H}_n} \| = {O}_{P^0_n} ( \sqrt{\frac{d_{\mathcal
{H}_n}}{N_n}} )$, as required.\vspace*{-2pt}
\end{pf*}

\begin{pf*}{Proof of Theorem \protect\ref{thmmodselconsistency}}
We follow the primal-dual witness method of Wainwright \cite{W06} as described in
Section \ref{secmodselconsistency}.
Let $\widetilde{\theta}_n \in\mathbb{R}^{d_{\mathcal{H}_n}}$ be the
restricted group-lasso estimator \eqref{eqtildetheta} and define the
vector $\hat{\theta}_n \in\mathbb{R}^{I-1}$, whose $h_n$ block is
given by
\[
\hat{\theta}_{h_n} = \cases{
\widetilde{\theta}_{h_n}, & \quad$\mbox{if } h_n \in\mathcal{H}_n$,\vspace*{1pt}\cr
0, & \quad$\mbox{otherwise.}$}\vspace*{-2pt}
\]
Then, by construction, $\widehat{\theta}_{\mathcal{H}_n} =
\widetilde
{\theta}_n$. Set also
\[
\widehat{\mathbf{m}}_n = N_n \frac{\exp(\mathrm{U} \widehat
{\theta
}_n)}{\langle
\exp(\mathrm{U} \widehat{\theta}_n), \mathbf{1} \rangle} = N_n
\frac{\exp(\mathrm{U}_{\mathcal{H}_n} \widetilde{\theta
}_n)}{\langle\\exp(\mathrm{U}_{\mathcal
{H}_n} \widetilde{\theta}_n), \mathbf{1} \rangle}.\vspace*{-2pt}
\]
Next, consider the random vector $\widehat{\eta} \in\mathbb{R}^{I-1}
=\operatorname{vec}(\widehat{\eta}_{\mathcal{H}_n},\widehat{\eta
}_{\mathcal
{H}_n^c})$, where
\[
\widehat{\eta}_{\mathcal{H}_n} = \operatorname{vec} \biggl\{
\lambda_{h_n} \frac{\widehat{\theta}_{h_n}}{ \| \widehat{\theta}_{h_n}
\|}, h_n \in\mathcal{H}_n \biggr\}\vspace*{-2pt}
\]
and
%
%
%e29 ###
\begin{equation}\label{eqzeta}
\widehat{\eta}_{\mathcal{H}_n^c} = \operatorname{vec} \{ \lambda_{w_n}
\widehat
{z}_{w_n}, w_n \in\mathcal{H}_n^c \} ,\vspace*{-2pt}
%, \| \widehat{z}_{h_n} \| \leq1.
\end{equation}
%
%and the matrix
with the sub-vectors $\{ \widehat{z}_{w_n}, w_n \in\mathcal{H}_n^c \}$
to be chosen in an appropriate way as described below. Notice also that
$\widehat{\eta}_{\mathcal{H}_n}$ belongs to the subdifferential of
$\lambda_n \sum_{h, h \in\mathcal{H}_n} \lambda_h \| x_h \|$ evaluated
at $\widetilde{\theta}_n$.

%The vector $\widehat{\eta}_{\mathcal{H}_n}$ is an explicit function of
%the group-lasso estimator.
% The matrix $\Sigma_{\mathcal{H}_n}$ is the submatrix of $\Sigma_n$
%indexed by the elements of $\mathcal{H}_n$ and is also $N_n$ times the
%Fisher information matrix at $\theta_n^0$, seen now as vector in $

The pair $(\widehat{\theta}_{\mathcal{H}_n} , \widehat{\eta
}_{\mathcal
{H}_n})$ must satisfy the optimality conditions (\ref{eqoptim}) for
the blocks indexed by $h_n \in\mathcal{H}_n$. Using this conditions
along with a Taylor expansion of $\widehat{\mathbf{m}}_n$ around
$\mathbf{m}^0_n$, we obtain the expression
%
%
%e30 ###
\begin{equation}\label{eqtheta}
\widehat{\theta}_{\mathcal{H}_n} = \theta^0_{\mathcal{H}_n} + N_n
\Sigma
_{\mathcal{H}_n}^{-1} \biggl( \frac{1}{N_n} \mathrm{U}^{\top}_{\mathcal
{H}_n}(\mathbf{n}_n - \mathbf{m}^0_n) - \frac{1}{N_n} \mathrm
{U}^{\top
}_{\mathcal
{H}_n} R_n - \lambda_n \widehat{\eta}_{\mathcal{H}_n} \biggr).\vspace*{-2pt}
\end{equation}
Employing a similar strategy, we now consider the optimality conditions
(\ref{eqoptim}) for the remaining blocks indexed by $w_n \in\mathcal
{H}_n^c$ and solve for $\{ \widehat{z}_{w_n}, w_n \in\mathcal{H}_n^c
\}
$ in terms of $\widehat{\theta}_{\mathcal{H}_n} - \theta
^0_{\mathcal
{H}_n}$. Eventually, we are led to the expression
%
%
%e31 ###
\begin{eqnarray}\label{eqz}
\lambda_n \widehat{\eta}_{\mathcal{H}_n^c} &=& \frac{1}{N_n} \mathrm
{U}^{\top
}_{\mathcal{H}_n^c}(\mathbf{n}_n - \mathbf{m}^0_n) - \frac{1}{N_n}
\mathrm{U}^{\top
}_{\mathcal{H}_n^c} R_n
\nonumber
\\[-10pt]
\\[-10pt]
\nonumber
&&{} - \mathrm{W}_n \biggl( \frac{1}{N_n} \mathrm
{U}^{\top
}_{\mathcal{H}_n}(\mathbf{n}_n - \mathbf{m}^0_n) - \frac{1}{N_n}
\mathrm{U}^{\top
}_{\mathcal{H}_n} R_n - \lambda_n \widehat{\eta}_{\mathcal{H}_n}
\biggr),\vspace*{-2pt}\vadjust{\goodbreak} %\lambda_n \frac{\widehat{\theta}^n_{H}}{\| \widehat{\theta}^n_{H}
\end{eqnarray}
where $\| R_n \| = {o}_{P_n^0}(\| \widehat{\theta}_{\mathcal{H}_n} -
\theta^0_{\mathcal{H}_n} \|)$, and
\[
\mathrm{W}_n = \mathrm{U}^{\top}_{\mathcal{H}_n^c} \biggl( D_{\mathbf
{m}^0_n} - \frac
{ \mathbf{m}^0_n (\mathbf{m}^0_n)^{\top}}{N_n}\biggr) \mathrm
{U}_{\mathcal{H}_n}
\Sigma_{\mathcal{H}_n}^{-1}. %( \m{U}^{\top}_{\mathcal{H}_n}
%( D_{\mathbf{m}^0_n} - \frac{\mathbf{m}^0_n (\mathbf{m}^0_n)^{
%) \m{U}_{\mathcal{H}_n} )^{-1}.
\]
Notice that, by Lemma \ref{lemestimationconsistecy}, $\| R_n \| =
{o}_{P^n_0}(1)$, so the remainder term in the above Taylor expansion is
negligible.

We then rely on equations (\ref{eqtheta}) and (\ref{eqz}) to show
that the assumed conditions are sufficient to guarantee that
%
%
%e32 ###
\begin{equation}\label{eqlimtheta}
\lim_n \mathbb{P} ( \| \widehat{\theta}_{h_n}\| > 0, \forall h_n
\in\mathcal{H}_n ) = 1
\end{equation}
and
%
%
%e33 ###
\begin{equation}\label{eqlimz}
\lim_n \mathbb{P} \Bigl( \max_{w_n \in\mathcal{H}_n^c}\| \widehat
{z}_{w_n} \| \leq1 \Bigr) =1.
\end{equation}

Because (\ref{eqlimtheta}) is equivalent to \eqref{eqOn1} and (\ref
{eqlimz}) is equivalent to \eqref{eqOn2}, model selection will follow.

%Using the pair $(\widehat{\theta}_{\mathcal{H}_n} , \widehat{\eta}_{
%around ${
%indexed by $h_n \in\mathcal{H}_n$ we obatine
%Using the optimality conditions (\ref{eqoptim}) and applying a Taylor
%expansion of $\widehat{\mathbf{m}}_n$ around $\mathbf{m}^0_n$, it can
%be
%verified that $\mathcal{O}_n$ holds if and only if both the equations
%Under the event $\mathcal{O}_n$, by applying a Taylor expansion of $
%optimality
%conditions (\ref{eqoptim}) as

%( \m{U}^{\top}_{\mathcal{H}_n} ( D_{\mathbf{m}^0_n} - \frac{{
%)

We will deal with equations (\ref{eqlimtheta}) and (\ref{eqlimz})
separately.

\textit{Proof of equation (\ref{eqlimtheta})}.
%We will show $\mathbb{P} ( \| \widehat{\theta}_{h_n}\| > 0,
It is enough to show that
%
%
%e34 ###
\begin{equation}\label{eqfirst}
%( D_{\mathbf{m}^0_n} - \frac{\mathbf{m}^0_n (\mathbf{m}^0_n)^{
%) \m{U}_{\mathcal{H}_n} )^{-1} ( \frac{1}{N_n} \m{U}^{
%) \| \leq\alpha_n ) \rightarrow1,
\mathbb{P} \biggl( \biggl\| N_n \Sigma_{\mathcal{H}_n}^{-1} \biggl( \frac
{1}{N_n} \mathrm{U}^{\top}_{\mathcal{H}_n}(\mathbf{n}_n - \mathbf
{m}^0_n) -
\frac
{1}{N_n} \mathrm{U}^{\top}_{\mathcal{H}_n} R_n - \lambda_n \widehat
{\eta
}_{\mathcal{H}_n} \biggr) \biggr\| \leq\alpha_n \biggr) \rightarrow1,
\end{equation}
where $\alpha_n = \min_{h_n \in\mathcal{H}_n} \| \theta^0_{h_n} \|$.
In fact, the former condition implies that the $h_n$-block of the
vector inside the norm sign in the previous display
%N_n ( \m{U}^{\top}_{H} ( D_\mathbf{m} - \frac{\mathbf{m}\mathbf{m}^{
%-
is less than $\| \theta^0_{h_n} \|$, $\forall h_n \in\mathcal{H}_n$,
which, by the triangle inequality, will produce the desired result.

%( \frac{1}{N_n} \m{U}^{\top}_{h_n}(\mathbf{n} - \mathbf{m}) -
%) \| < \| \theta^0_{h_n} \|, ,\| \theta^0_{h_n} \|,
%) \rightarrow1,
%which, by the triangle inequality, will produce the desired result.

%This involves, more generally, looking at the matrix $( \m{U}^{
%complicated.

%D_{\mathbf{m}^0_n} - \frac{\mathbf{m}^0_n(\mathbf{m}^0_n)^{\top}}{N_n}
%)

First, we consider the term
\[
\Sigma_{\mathcal{H}_n}^{-1} \mathrm{U}^{\top}_{\mathcal
{H}_n}(\mathbf{n}_n -
\mathbf{m}^0_n).
\]
The vector $\mathrm{U}^{\top}_{\mathcal{H}_n}(\mathbf{n}_n -
\mathbf
{m}^0_n)$ has
mean zero and covariance matrix $\Sigma_{\mathcal{H}_n}$. Furthermore,
because of [NC.2], letting $\gamma^{\min}_{\mathcal{H}_n} = \lambda
_{\min} ( \Sigma_{\mathcal{H}_n} ) $, we have
%
%
%e35 ###
\begin{equation}\label{eqeig}
\gamma^{\min}_{\mathcal{H}_n} \frac{1}{N_n} \geq D_{\min} > 0\qquad
\mbox{for all } n.
\end{equation}
%
%We want to show that $\| I_1 \|$ is ${O}_P(r_n)$ for some $r_n$ to be
%determined.
%Then, (\ref{eqft}) has mean zero and covariance matrix letting $X_n =
Combining these observations, and using the formula for the expected
value of a quadratic form, we arrive at
\[
\mathbb{E} \| \Sigma_{\mathcal{H}_n}^{-1} \mathrm{U}^{\top
}_{\mathcal
{H}_n}(\mathbf{n}_n - \mathbf{m}^0_n) \|^2_2 %= \mathbb{E} X^{\top}_n
= \operatorname{tr} \Sigma^{-1}_{\mathcal{H}_n} \leq\frac{d_{\mathcal
{H}_n}}{\gamma^{\min}_{\mathcal{H}_n}} \leq\frac{d_{\mathcal
{H}_n}}{D_{\min} N_n}, %\leq\frac{d_{\mathcal{H}_n}}{\lambda_n^{\min}} ,
\]
%
%= \m{tr} ( \m{U}^{\top}_{H} ( D_\mathbf{m} - \frac{\mathbf{m}
where $d_{\mathcal{H}_n} = \sum_{h \in\mathcal{H}_n} d_h$.
% and $\gamma_n^{\min}$ is the smallest eigenvalue of $\Sigma_{
%, which is not smaller than the smallest eigenvalue $\lambda_n^{\min}$
%of the Fisher information matrix.
Then, Chebyshev inequality implies
%
%
%e36 ###
\begin{equation}\label{eqfirst1}
\| \Sigma_{\mathcal{H}_n}^{-1} \mathrm{U}^{\top}_{\mathcal
{H}_n}(\mathbf
{n}_n - \mathbf{m}^0_n) \| = {O}_{P_n^0}\Biggl( \sqrt{\frac{d_{\mathcal
{H}_n}}{N_n}} \Biggr).
\end{equation}
%
%information matrix to the one of the Fisher information matrix compute
%using the whole $\m{U}$.

Next, using (\ref{eqeig}) for the operator norm of of $\Sigma
_{\mathcal
{H}_n}$, we get the upper bound
%
%
%e37 ###
\begin{equation}\label{eqfirst2}
\| N_n \Sigma_{\mathcal{H}_n}^{-1} \lambda_n \widehat{\eta
}_{\mathcal{H}_n} \| \leq\frac{1}{ D_{\min}} \lambda_n \sqrt{\sum
_{h_n \in H} \lambda^2_{h_n}},
\end{equation}
%
%is such that
%so that we need
which, for $\lambda_{h_n} = \sqrt{d_{h_n}}$, simplifies to $\frac{1}{
D_{\min}} \lambda_n \sqrt{d_{\mathcal{H}_n}}$.

Finally, the norm of
\[
\Sigma_{\mathcal{H}_n}^{-1} \mathrm{U}^{\top}_{\mathcal{H}_n} R_n
\]
is no larger than
%
%
%e38 ###
\begin{equation}\label{eqfirst3}
\frac{1}{\sqrt{D_{\min} N_n}} \sqrt{d_{\mathcal{H}_n}} {o}(\|
\widehat{\theta}_{\mathcal{H}_n} - \theta^0_{\mathcal{H}_n} \|
) = {o}_{P_n^0} \biggl( \frac{d_{\mathcal{H}_n}}{N_n} \biggr),
\end{equation}
because $\| \mathrm{U}_{\mathcal{H}_n} \| \leq\sqrt{d_{\mathcal{H}_n}}$.

Using equations (\ref{eqfirst1}), (\ref{eqfirst2}) and (\ref
{eqfirst3}), condition (\ref{eqfirst}) is satisfied if MSC.1 holds.\\
%which, for $\lambda_{h_n} = \sqrt{d_{h_n}}$, simplifies to $\frac{

\textit{Proof of equation (\ref{eqlimz}).}
%We need to show %that the blocks of $\widehat{z}_{H_n^c}$ have all
%norms no greater than $1$ with increasingly (in $n$) large
%probability, i.e. that
%We will show the stronger result $\lim_n Pr \{ \| \widehat{z}_{H^c}
%We need this assumption:
%is
%similar to the assumptions made on the design matrix for Gaussian
%ensemble and implies $\| \m{W}_n \|_{2} < 1 - \varepsilon$.\\=
%Actually, we may only need this assumption:
%and
%each $w_n \in H_n^c$,%$h \in H$ and
%(\mathbf{m}^0_n)^{\top}}{N_n}) \m{U}_{H_n} ( \m{U}^{
%%\frac{d_h}{d_H}
%or the more stringent one
%(\mathbf{m}^0_n)^{\top}}{N_n}) \m{U}_{H_n} ( \m{U}^{
In equation (\ref{eqz})
write $\widehat{\eta}_{\mathcal{H}_n^c} = \Lambda_{\mathcal{H}_n^c}
\widehat{z}_{\mathcal{H}_n^c}$, %and $\widehat{\eta}^n_{\mathcal{H}_n}
%=
where $\Lambda_{\mathcal{H}_n^c}$ is a $d_{\mathcal
{H}_n^c}$-dimensional diagonal matrix whose diagonal is $\operatorname{vec}\{
\mathbf{1}_{w_n}\lambda_{w_n}, w_n \in\mathcal{H}^c_n\}$, with
$\mathbf{1}_{h_n}$
denoting the $d_{h_n}$-dimensional vector with entries all equal to
$1$. Then, (\ref{eqz}) becomes
% and, by the optimality conditions, $\| \widehat{z}_{h} \| \leq1 $,
%$h \in H^c$, to obtain
%
\begin{eqnarray*}
\widehat{z}_{\mathcal{H}_n^c} &=& \frac{1}{N_n \lambda_n} \Lambda
_{\mathcal{H}_n^c}^{-1} \mathrm{U}^{\top}_{\mathcal{H}_n^c}(\mathbf
{n}_n -
\mathbf{m}^0_n) - \frac{1}{N_n \lambda_n} \Lambda_{\mathcal
{H}_n^c}^{-1} \mathrm{U}^{\top}_{\mathcal{H}_n^c}R_n\\
&&{} - \Lambda
_{\mathcal{H}_n^c}^{-1} \mathrm{W}_n \biggl( \frac{1}{N_n \lambda_n}
\mathrm{U}^{\top}_{\mathcal{H}_n}(\mathbf
{n}_n - \mathbf{m}^0_n) - \frac{1}{N_n \lambda_n} \mathrm{U}^{\top
}_{\mathcal
{H}_n} R_n - \widehat{\eta}_{\mathcal{H}_n} \biggr).
\end{eqnarray*}

For any $w_n \in\mathcal{H}_n^c$, consider the corresponding block in
the vector $\Lambda_{\mathcal{H}_n^c}^{-1} \mathrm{W}_n \widehat
{\eta
}_{\mathcal{H}_n}$, that is, the vector
%
%
%e39 ###
\begin{equation}\label{eqlalala}
\frac{1}{\lambda_{w_n}} \mathrm{U}^{\top}_{w_n} \Sigma^0_n \mathrm
{U}_{\mathcal
{H}_n} ( \mathrm{U}^{\top}_{\mathcal{H}_n}\Sigma^0_n \mathrm
{U}_{\mathcal
{H}_n} )^{-1} \widehat{\eta}_{\mathcal{H}_n}.%\frac{\widehat{
\end{equation}
Because of assumption [MSC.2], the Euclidian norm of (\ref{eqlalala}),
for any choice of $w_n \in\mathcal{H}_n^c$, is bounded by
\[
\frac{(1- \varepsilon)}{|\mathcal{H}_n^c|} \frac{\sum_{h_n \in
\mathcal
{H}_n} \lambda_{h_n}}{\min_{w_n \in\mathcal{H}_n^c}\lambda_{w_n}},
\]
which, in turn, is smaller that
\[
(1- \varepsilon) \frac{|\mathcal{H}_n| \max_{h_n \in\mathcal
{H}_n} \lambda
_{h_n}}{|\mathcal{H}_n^c| \min_{w_n \in\mathcal{H}_n^c}\lambda_{w_n}}.
\]
Then, under MSC.3 (\ref{eqlalala}) will be eventually less than
$(1-\varepsilon)$, uniformly over $w_n \in\mathcal{H}_n^c$.

%Alternatively, under (\ref{eqparort2}), the bound becomes
%Notice that, if $\lambda_h =1$ for each $h$, then the above bound
%would just be $(1-\varepsilon)$.

%Then, the sparsity assumption on the penalties $\lambda_h$ when (
%or, if (\ref{eqparort2}) is in force instead,
%Notice that if $\lambda_{h_n} = \sqrt{d_{h_n}}$ the last display is
%always satisfied for each $n$.

Next, for $w_n \in\mathcal{H}_n^c$, we consider the vector
%
%
%e40 ###
\begin{equation}\label{eqpoipoipoi}
\frac{1}{N_n \lambda_n \lambda_{w_n}} [ \mathrm{U}^{\top
}_{w_n}(\mathbf{n}_n
- \mathbf{m}^0_n) - \mathrm{W}_{w_n} \mathrm{U}^{\top}_{\mathcal
{H}_n}(\mathbf
{n}_n - \mathbf{m}^0_n) ].
\end{equation}
The covariance matrix of the term inside the parenthesis is
%
%
%e41 ###
\begin{equation}\label{eqlargeigen}
\mathrm{U}_{w_n}^{\top} (\Sigma^0_n)^{1/2} [ I_{d_{\mathcal{H}_n}} -
(\Sigma^0_n)^{1/2} \mathrm{U}_{\mathcal{H}_n} ( \mathrm
{U}_{\mathcal
{H}_n}^{\top} \Sigma^0_n \mathrm{U}_{\mathcal{H}_n} )^{-1} \mathrm
{U}_{\mathcal{H}_n}^{\top} (\Sigma^0_n)^{1/2} ] (\Sigma
^0_n)^{1/2} \mathrm{U}_{w_n},
\end{equation}
where
\[
\Sigma^0_n = D_{\mathbf{m}^0_n} - \frac{\mathbf{m}^0_n (\mathbf
{m}^0_n)^{\top}}{N_n}.
\]
%
%Since the largest eigenvalue of the matrix in (\ref{eqlargeigen}) is
%smaller than $N_n l^{\max}_n$, by Chebyshev inequality, the term (
%Under the condition [MSC.4], uniformly over $w_n \in\mathcal{H}_n^c$,
%the expression (\ref{eqpoipoipoi}) is ${o}_{P^0_n}(1)$.

%OLD
%is no greater than
% \frac{1}{\lambda_n \lambda_{w_n} N_n} [ \| \m{U}_{w_n}^{
% which is of order
%) ],
%by Chebyshev inequality.

%For every $\delta> 0$, using the fact that $\| \m{U}^{\top}_{h_n} ({
%Bretagnolle-Huber-Carol inequality \citep[see, e.g.,][Proposition
%A.6.6]{VVW:98} we obtain
%and
%)^2 \frac{N_n}{d_{w_n}} + I_n \log2 \}.%3 \exp\{ -
%- \mbf{\pi}^0_n) \|_{2} \leq\sqrt{d_{h_n}} \| (\mathbf{p}_n - \mbf{

%to $0$.}

%Therefore, under the condition MSC.4, the norm of vector (
%) \rightarrow\infty

%( \lambda_n \lambda_{w_n})^2 \frac{N_n}{d_{H_n}}
%and
%( \lambda_n \lambda_{w_n})^2 \frac{N_n}{d_{w_n}}
%convergence in probability to zero is granted. Notice that, for the
%choice of penalty terms $\lambda_{h_n} = \sqrt{d_{h_n}}$, the previous
%condition simplifies to
% \lambda_n^2 N_n \rightarrow\infty

Since the largest eigenvalue of the matrix in (\ref{eqlargeigen}) is
smaller than the largest eigenvalue of the covariance matrix of
$\mathrm{U}_{w_n}^\top(\mathbf{n}_n - \mathbf{m}^0_n)$, by Chebyshev's
inequality it
is enough to show that the $\ell_2$ norm of
\[
\frac{1}{N_n \lambda_n \lambda_{w_n}} \mathrm{U}^{\top
}_{w_n}(\mathbf{n}_n -
\mathbf{m}^0_n)
\]
vanishes in order to conclude that (\ref{eqpoipoipoi}) has vanishing
$\ell_2$ norm as well.
To this end, notice that
\begin{eqnarray*}
\frac{1}{N_n \lambda_n \lambda_{w_n}} \| \mathrm{U}^{\top
}_{w_n}(\mathbf
{n}_n -
\mathbf{m}^0_n) \| &\leq&\frac{\sqrt{d_{w_n}}}{N_n \lambda_n \lambda_{w_n}}
\| \mathrm{U}^{\top}_{w_n}(\mathbf{n}_n - \mathbf{m}^0_n) \|_\infty
\\
&\leq&\frac
{\sqrt{d_{w_n}}}{N_n \lambda_n \lambda_{w_n}} \| \mathrm{U}^{\top
}_{\mathcal
{H}^c_n}(\mathbf{n}_n - \mathbf{m}^0_n) \|_\infty.
\end{eqnarray*}
Next, write
\[
\mathrm{U}^{\top}_{\mathcal{H}^c_n} \frac{(\mathbf{n}_n - \mathbf
{m}^0_n)}{N_n} =
\sum_{j_n =1}^{N_n} \mathrm{U}^{\top}_{\mathcal{H}^c_n} \frac
{(X_{j_n} - \bolds{\pi}^0_n)}{N_n},
\]
where the vectors $X_{j_n}, 1 \leq j_n \leq N_n$ are i.i.d.
$\operatorname{Multinomial}(1,\bolds{\pi}^0_n)$. Since the entries of
$\mathrm{U}_{\mathcal
{H}^c_n}$ are all $-1$, $0$ or $1$, by Bernstein's inequality followed
by a union bound, we get
\[
\mathbb{P} \biggl\{ \biggl\| \frac{1}{N_n} \mathrm{U}^{\top}_{\mathcal
{H}^c_n}(\mathbf{n}_n - \mathbf{m}^0_n) \biggr\|_\infty> c \frac{\lambda_n
\lambda_{w_n}}{\sqrt{d_{w_n}}} \biggr\} \leq2 \exp\biggl\{ - \frac{N_n
c^2 \lambda_n^2 {\lambda_{w_n}^2}/{d_{w_n}}}{1/8 + ({2}/{3}) c
{\lambda_n \lambda_{w_n}}/{\sqrt{d_{w_n}}}} + \log d_{\mathcal
{H}^c_n}\biggr\},
\]
which vanishes under [MSC.4].
As for the terms involving $R_n$, following the arguments used above,
it is easy to see that they both converge in probability to $0$, so
that (\ref{eqlimz}) holds true.
\end{pf*}

\begin{pf*}{Proof of Theorem \protect\ref{thmclt}}
All the claims in the proof are made on the event $\mathcal{O}_n$.
Because the norm consistency assumptions are in force, $\mathcal{O}_n$
occurs in probability and, therefore, our claims hold true within a set
or probability converging to $1$.
In particular, $\| \widehat{\theta}_{\mathcal{H}_n} - \theta
^0_{\mathcal
{H}_n}\| = {O}_{P_n^0} ( \sqrt{\frac{d_{\mathcal{H}_n}}{N_n}}
) (1+{o}_{P_n^0}(1)) = {O}_{P_n^0} ( \sqrt{\frac{d_{\mathcal
{H}_n}}{N_n}} )$.
Reorganize equation (\ref{eqtheta}) as
%
%
%e42 ###
\begin{equation}\label{eqcltrepresentation}
\Sigma_{\mathcal{H}_n}^{1/2} ( \widehat{\theta}_{\mathcal{H}_n} -
\theta^0_{\mathcal{H}_n} )  =  \Sigma_{\mathcal{H}_n}^{-1/2}
\mathrm{U}^{\top}_{\mathcal{H}_n}(\mathbf{n}_n - \mathbf{m}^0_n) -
\Sigma
_{\mathcal
{H}_n}^{-1/2} N_n \lambda_n \widehat{\eta}_{\mathcal{H}_n} - \Sigma
_{\mathcal{H}_n}^{-1/2} \mathrm{U}^{\top}_{\mathcal{H}_n} R_n.
%&\equiv& Z_n+\{R^1_n+R^2_n\}.
\end{equation}

By similar arguments used in the proof of Theorem \ref
{thmmodselconsistency}, the term
\[
\Sigma_{\mathcal{H}_n}^{-1/2} \mathrm{U}^{\top}_{\mathcal{H}_n} R_n
\]
is of order
\[
\sqrt{\frac{d_{\mathcal{H}_n}}{N_n}} {o}_{P_n^0}( \| \widehat{\theta
}_n - \theta^0_n \| ) = {o}_{P_n^0}\biggl( \frac{d_{\mathcal
{H}_n}}{N_n}\biggr),
\]
and therefore converges in probability to $0$.

%its norm is bounded by
%because $\|\hat\eta_{\mathcal{H}_n}\| = \sqrt{\sum_{h_n \in
%eigenvalue of $\Sigma_{\mathcal{H}_n}$. In order for (
%eventually guarantee the same type of asymptotic optimality as
%sequences of MLEs, withotu actually knowing the true underlying model,
%in practice, it may be too strong of a penalty.
%We commence with $R_n^1$ and show that this term is $o(1)$ as $n
%vector comprised by $h_n$-blocks of the form $\lambda_{h_n} \hat
%by assumption \textit{(i)}, thus showing $R^1_n$ is asymptotic
%negligible. \ale{The previous display becomes $\lambda_n \sqrt{N_n}
%using Chebyshev inequality,
%I_n^{1/2} N_n^{1/2}\rightarrow0 , \mbox{for all}
%which proves that $R_1$ tends in probability to zero.$d_{H_n}^{1/2}$
%Alternatively, a much weaker condition on the rate of decay of $
%adding a bias term to both the asymptotic mean and variance of $
As for $\Sigma_{\mathcal{H}_n}^{-1/2} N_n \lambda_n \widehat{\eta
}_{\mathcal{H}_n}$,
notice that, on $\mathcal{O}_n$, the vector $\widehat{\eta}^0_{H_n}$ is
a differentiable function of $\widehat{\theta}_{\mathcal{H}_n} \in
\mathbb{R}^{d_{\mathcal{H}_n}}$.
Then, using a Taylor expansion around $\theta^0_{\mathcal{H}_n}$,
%
%
%e43 ###
\begin{equation}\label{eqtaylorpenalty}
\Sigma^{-1/2}_{\mathcal{H}_n} N_n \lambda_n \widehat{\eta
}_{\mathcal
{H}_n} = \Sigma^{-1/2}_{\mathcal{H}_n} N_n \lambda_n \Biggl( \eta
^0_{\mathcal{H}_n} + \mathrm{J}^0_{\mathcal{H}_n} ( \widehat{\theta
}_{\mathcal{H}_n} - \theta^0_{\mathcal{H}_n} ) + {o}_{P_n^0}\Biggl(
\sqrt{\frac{d_{\mathcal{H}_n}}{N_n}} \Biggr) \Biggr). %{o}_P( \sqrt{
\end{equation}
%
%where $\m{J}^0_{\mathcal{H}_n}$ is a $d_{\mathcal{H}_n} \times d_{
%$d_{h_n} \times d_{h_n}$ matrix
%with $\m{I}_{d_{h_n}}$ denoting the $d_{h_n}$-dimensional identity
%matrix.
The remainder term in equation (\ref{eqtaylorpenalty}) is of order
\[
\lambda_n {o}_{P_n^0}(\sqrt{d_{\mathcal{H}_n}}),% \sqrt{I_n}),
\]
which become negligible for $\lambda_n = {O} ( \frac{1}{\sqrt
{d_{\mathcal{H}_n}}} )$ (obviously, $\lambda_n = {O} ( 1/\sqrt
{I_n} )$ will do). %However, we now add an extra bias to both the
%mean and the variance of $\widehat{\theta_n^0}$, because we now we
%will establish normality of this sequence instead
Then using~(\ref{eqcltrepresentation}), we obtain %can see that $
%probability, to
%
%
%e44 ###
\begin{equation}\label{eqZn}
%( \Sigma^{1/2}_{H_n} + \Sigma^{-1/2}_{H_n} N_n \lambda_n
%J^0_{H_n} ) ( \widehat{\theta}_n - \theta^0_n ) +
\Sigma_{\mathcal{H}_n}^{-1/2} \mathrm{U}^{\top}_{\mathcal
{H}_n}(\mathbf{n}_n -
\mathbf{m}^0_n) = \Sigma_{\mathcal{H}_n}^{-1/2}\bigl ( ( \Sigma
_{\mathcal{H}_n} + N_n \lambda_n J^0_{\mathcal{H}_n} ) ( \widehat
{\theta}_n - \theta^0_n ) + N_n \lambda_n \eta^0_{\mathcal{H}_n}
\bigr) + {o}_{P^0_n}(1).
\end{equation}
%
%of the following
%N_n \lambda_n J^0_{H_n} )$.
%Bringing the $\Sigma^{-1/2}_{H_n} N_n \lambda_n N_n \eta^0_{H_n}$ on
%the other side and notice that the reminder term is
%Another option is to deal directly with the term $\Sigma^{-1/2}_{H_n}
%N_n \lambda_n \widehat{\eta}_{H_n}$ and notice that its norm is
%which, for $\lambda_{h_n} = \sqrt{d_{h_n}}$ gives $\lambda_n \sqrt{N_n
%d_{H_n}} = o(1)$. We could even continue with a second order
%expansion, noting that In case this is useful, for $x \neq0$, $\nabla
%( \m{I} - \frac{x}{\Vert x\Vert _2} ( \frac{x}{\Vert x\Vert _2} )^{
%require $\lambda_n = O( \frac{N_n}{\sqrt{I_n}})$.}
Thus, we only need to consider the term $\Sigma_{\mathcal{H}_n}^{-1/2}
\mathrm{U}^{\top}_{\mathcal{H}_n}(\mathbf{n}_n - \mathbf{m}^0_n)$.
For $1
\leq j_n
\leq N_n$, let
\[
Y_{j_n} = \frac{1}{\sqrt{N_n}} \mathrm{F}_n^{-1/2} \mathrm{U}^{\top
}_{\mathcal
{H}_n} (X_{j_n} - \bolds{\pi}^0_n),
\]
where the variables $X_{j_n}$ are \mbox{i.i.d.} Multinomials with size
$1$ and probability vector $\bolds{\pi}^0_n$.
Then,
\[
\Sigma_{\mathcal{H}_n}^{-1/2} \mathrm{U}^{\top}_{\mathcal
{H}_n}(\mathbf{n}_n -
\mathbf{m}^0_n) = \sum_{j_n} Y_{j_n},
\]
where $\mathbb{E}Y_{j_n} = 0$, $\operatorname{Cov} Y_{j_n} = \frac
{1}{N_n}\mathrm{I}_{d_{\mathcal{H}_n}} $ and
$\sum_{j_n} \operatorname{cov}(Y_{j_n}) = \mathrm{I}_{d_{\mathcal{H}_n}}$.

To show the result in part 1 it is sufficient to show that, under
assumption [CLT.LF], the multivariate Lindberg--Feller conditions hold, namely
\[
\sum_{j_n} \mathbb{E}_{\theta^0_{\mathcal{H}_n}} \| \mathrm{G}_n
Y_{j_n} \|
^2 I_{\{ \| \mathrm{G}_n Y_{j_n} \| \geq\varepsilon\}} \rightarrow0
\]
as $n \rightarrow\infty$, where $I_{\{ \cdot\}}$ denotes the indicator
function. The proof is quite standard (see also the proof of
Theorem 2 in Fan and Peng \cite{fan-2004-32}) and we only sketch it.

% We omit this proof because it is almost identical to the proof we
%produce below for the [CLT.LF] conditions, the only difference being a
%rate $O(N_n^{-1})$ in equation (\ref{eqche}). See also the proof of
%Theorem 2 in \cite{fan-2004-32}.

Because the vectors $Y_{j_n}$'s are identically distributed, and
invoking the Cauchy--Schwarz inequality, it is sufficient to show that
%
%
%e45 ###
\begin{equation}\label{eqcltlf}
N_n \bigl( \mathbb{E} \| \mathrm{G}_n Y_{j_n} \|^4 \mathbb{P}(\| \mathrm{G}_n
Y_{j_n} \| \geq\varepsilon) \bigr)^{1/2}\rightarrow0.
\end{equation}

By Chebychev inequality, for fixed $\varepsilon> 0$,
\[
\mathbb{P}(\| \mathrm{G}_n Y_{j_n} \| \geq\varepsilon) \leq\frac
{\operatorname{tr}(\mathrm{G}_n\mathrm{G}_n^\top)}{\varepsilon^2 N_n} =
{O} \biggl( \frac{1}{N_n} \biggr),
\]
where $\operatorname{tr}(\mathrm{G}_n\mathrm{G}_n^\top) = {O}(1)$ because
of \eqref{eqGn}.
Similarly, using the fact that the minimal eigenvalue of $\mathrm
{F}_n$ is
bounded away from zero,
\[
\mathbb{E} \| \mathrm{G}_n Y_{j_n} \|^4 \leq {O}\biggl(\frac{1}{N_n^2}\biggr)
\mathbb{E} \| \mathrm{U}^{\top}_{\mathcal{H}_n}(X_{j_n} - \bolds
{\pi}^0_n) \|
^4 = {O}\biggl( \frac{d_{\mathcal{H}_n}^{2}}{N_n^{2}}\biggr),
\]
where in the last step we use the fact that the entries of $\mathrm
{U}^{\top
}_{\mathcal{H}_n}(X_{j_n} - \bolds{\pi}^0_n)$ are bounded, uniformly
over $n$.
%where the norm in last expectation is the Frobenius norm and we used
%the fact that $x^{\top} x \leq(\sum_{i,j} |x_i x_j|^2)^{1/2}$. Since,
%the above term is of order $O( \frac{d_{
%Similarly, using the formulas for the mean and variance of a quadratic
%form, for each $j_n$,
%)
Combining the last two displays, the left-hand side of (\ref
{eqcltlf}) is of order
\[
N_n {O} \Biggl( \sqrt{\frac{d_{\mathcal{H}_n}}{N_n}} \frac{1}{N_n} \Biggr)
= {O} \biggl( \frac{d_{\mathcal{H}_n}}{N_n^{1/2}}\biggr),
\]
which, in virtue of assumption [CLT.LF], vanishes, as desired.

Next, we prove the result of part 1 under both [CLT.Ma] and [CLT.Mb].
%First we show this under the assumption [CLT.LF], which allows for a
%classic Lindeberg-Feller argument. %Assume that
%$\max_{i \in\mathcal{I}_n} \pi_i^0 \leq D_1 (N_n d_{
%%$\max_{i \in\mathcal{I}_n} \pi_i^0\leq D_1 I_n^{-1}$, for some
%positive constant $D_1$;
%%\item[CLT.2] %(ii)
%%$d_{\mathcal{H}_n}=o( I_n )$;
%$\lambda_n = {O} ( \frac{1}{\sqrt{d_{\mathcal{H}_n}}} )$;
%$d_{\mathcal{H}_n} = o ( \sqrt{N_n} )$.
We relax the assumption [CLT.LF] by allowing the dimension of the
parameter space to grow fatser. To this end, we derive
multi-dimensional analogs of Lemmas 2.1 and 2.2 and Theorem 2.1 in
Morris~\cite{MR0370871}. In particular, our proof follows closely the proof of
Morris \cite{MR0370871}, Lemma 2.2. %The idea is basically the same.
We first obtain joint limit law by using Lemma~\ref{lmlemma21}, and
then establish the conditional limit law by using a multi-dimensional
version of condition $(2.9)$ in~Morris \cite{MR0370871}. Note that the result
in~Steck \cite{MR0091552} about conditional limit laws is actually a
multi-dimensional one, but somehow was formulated in Morris \cite{MR0370871},
Theorem~2.1, as one-dimensional. The conditional law we are
interested is the distribution of $Z_n$, defined below in (\ref{eqsteck}).

Let $\gamma_n = N_n^{-1}\mathrm{G}_n\Sigma_{\mathcal
{H}_n}^{-1/2}\mathrm{U}_{\mathcal{H}_n}^{\top}\mathbf{m}_n^0$,
and set $\mathrm{A}_n=\mathrm{G}_n\Sigma_{\mathcal
{H}_n}^{-1/2}\mathrm{U}_{\mathcal
{H}_n}^{\top}$.
Note that $\mathbf{m}^0_n=N_n \bolds{\pi}_n^0$, thus $\gamma
_n=\mathrm{A}_n
\bolds{\pi}_n^0$. Denote the $i$th column of $\mathrm{A}_n$ by
$a_i$, $i=1,\ldots,
I_n$. Then, the left-hand
side of (\ref{eqZn}), premultiplied by $\mathrm{G}_n$, can be written as
\[
Z_n = \sum_{i \in\mathcal{I}_n} f_i(n_i),
\]
where $f_i(n_i)=(a_i-\gamma_n)(n_i-m^0_i)$.
Let $\{X_i ; i=1,\ldots, I_n\}$ be independent Poisson random variables
with mean $m_i^0=N_n \bolds{\pi}_i^0$, so that $\mathbb
{E}f_i(X_i)=0$ and
$\sum_{i}\operatorname{cov}(f_i(X_i),X_i)=0$, by construction.
Next, define
%
%
%e47 ###
%e46 ###
\begin{eqnarray}
V_n &=& N_n^{-1/2}\sum_{i}(X_i-m_i^0), \\
U_n &=& \Xi_n^{-1/2}\sum_{i}f_i(X_i) ,
\end{eqnarray}
where $\Xi_n = \sum_{i} \operatorname{cov} (f_i(X_i))$. A simple calculation
though shows that $\Xi_n = \mathrm{G}_n\mathrm{G}_n^{\top}$,
a~square matrix of fixed dimensions $k\times k$.
The goal is to prove the asymptotic normality of $U_n$ given $\{V_n=0\}
$, and then
use the fact (underlying Morris' method) that
%
%
%e48 ###
\begin{equation}\label{eqsteck}
\mathcal{L}(\Xi_n^{-1/2} Z_n) =\mathcal{L}(U_n|V_n=0), %\mathcal{L}(
\end{equation}
where $\mathcal{L}$ stands for law.

The random variables $V_n$ have zero means and unit variances.
Furthermore, by the
same arguments used in the early parts of (Morris \cite{MR0370871}, Lemma 2.2),
assumption [CLT.Mb] guarantees that the \textit{uan} condition is
satisfied, so the sequence $V_n$ converge in distribution to a Gaussian
variable.
%Lienberg condition is satisfied.
%N_n^{-1} \sum_{i \in\mathcal{I}_n} \mathbb{E} ( | X_i - m^0_i
%|^2 I\{ | X_i - m^0_i | > \varepsilon\sqrt{N_n} \} ) \rightarrow
%0, \mbox{ as } n \rightarrow\infty
%By Cauchy Swartz the previos term is less than
%N_n^{-1} \sqrt{ \sum_{i \in\mathcal{I}_n} (m^0_i - 3 (m^0_i)^2) \exp
%}
Similarly, the random vector $U_n$ satisfies
$\mathbb{E}U_n=0$, $\operatorname{cov}(U_n)= \mathrm{I}_k$, the identity
matrix of dimensions $k\times k$, and, by construction, $\operatorname
{cov}(V_n,U_n)=0$. We argue below that $U_n$ satisfies the
multi-dimensional Lindeberg condition. By Lemma~(\ref{lmlemma21}),
this will imply the asymptotic normality of the joint limit law of $(V_n,U_n)$.

By Schwartz inequality, for any $\varepsilon>0$,
%
%
%e49 ###
\begin{equation}\label{eqlindschwarz}
\sum_i \mathbb{E}[ \|f_i(X_i) \|^2 ; \|f_i(X_i) \|>
\varepsilon] \leq\sum_i \bigl[\mathbb{E}\|f_i(X_i) \|^4
\mathbb{P}\bigl(\|f_i(X_i) \|> \varepsilon\bigr) \bigr]^{1/2} .
\end{equation}
We will show that, for each $\varepsilon>0$, the right-hand
side of~(\ref{eqlindschwarz}) tends to zero.
Recall that $f_i(X_i)=(a_i-\gamma_n)(X_i-m_i^0)$.
The length of $\gamma_n$ can be bounded as follows,
$\|\gamma_n\| \leq\|\mathrm{G}_n\|
\|\Sigma_{\mathcal{H}_n}^{-1/2}\mathrm{U}_{\mathcal{H}_n}^{\top}
\bolds{\pi }_n^0\|
\leq {O}(1) D_{\operatorname{min}} N_n^{-1/2} \| \mathrm{U}_{\mathcal
{H}_n}^{\top} \bolds{\pi}_n^0\|$.
Elements of $\mathrm{U}_{\mathcal{H}_n}^{\top} \bolds{\pi}_n^0$
are absolutely
bounded by a
constant $D_1$, thus $\|\gamma_n \| \leq D N_n^{-1/2} d_{\mathcal
{H}_n}^{1/2}$.
Similarly, $\|a_i\|=\|\mathrm{A}_n e_i\|\leq D N_n^{-1/2} d_{\mathcal
{H}_n}^{1/2}$,
where $e_i$ is the standard unit vector in $\mathbb{R}^{I_n}$ with $i$th
coordinate equal to $1$. Adding up,
$\|a_i-\gamma_n \|={O} (N_n^{-1/2}d_{\mathcal{H}_n}^{1/2} )$,
which tends to zero by assumption [CLT.Ma].

Next, we use the following large deviation result for Poisson random
variables, due to Bobkov and Ledoux \cite{MR1636948} and based on a modified logarithmic
Sobolev inequality:
\begin{theorem}\label{thmsobolev}
Let $X$ be a Poisson random variable with parameter $\lambda$.
Then, for every $h:\mathbb{N}\rightarrow\mathbb{R}$, with
$\sup_{x\in\mathbb{N}}|h(x+1)-h(x)|\leq1$,
%
%
%e50 ###
\begin{equation}\label{eqsobolev}
\mathbb{P}\bigl( h(X)-\mathbb{E}h(X) \geq b \bigr) \leq
\exp\biggl\{-\frac{b}{4}\log\biggl(1+\frac{b}{2\lambda}\biggr)\biggr\}
,
\end{equation}
for all $b\geq0$.
\end{theorem}

Then, using Theorem \ref{eqsobolev}, for some constant $D$,
%
%
%e51 ###
\begin{eqnarray}\label{eqsobolevpoisson}
&&\mathbb{P}(X_i-m_i^0\geq\varepsilon\|a_i-\gamma_n \|^{-1}
)\nonumber\\
&&\quad\leq \exp\biggl\{-\frac{\varepsilon}{4}\|a_i-\gamma_n \|^{-1} \log\biggl(
1+\frac{1}{2}\frac{\varepsilon}{m_i^0\|a_i-\gamma_n\|} \biggr)\biggr\}
\nonumber
\\[-8pt]
\\[-8pt]
\nonumber
&&\quad\leq \exp\biggl\{-\varepsilon D N_n^{1/2}d_{\mathcal{H}_n}^{-1/2} \log
\biggl(1+\varepsilon D \frac{1}{\sqrt{N_n d_{\mathcal{H}_n}} \max_i
\bolds{\pi }^0_i} \biggr)\biggr\} \\
&&\quad =  \exp\bigl( -{O} \bigl( \sqrt{N_n /d_{\mathcal{H}_n}} \bigr) \bigr)
,\nonumber
%&\leq& \exp\{-\varepsilon D N_n^{1/2}d_{\mathcal{H}_n}^{-1/2} \log
%(1+\varepsilon D \frac{I_n}{N_n^{1/2}d_{\mathcal{H}_n}^{1/2}} )
%&\leq& \OO{e^{-I_n/d_{\mathcal{H}_n}}} \longrightarrow0
\end{eqnarray}
as $n\rightarrow\infty$.
%The last inequality follows by noting that $\log(1+x)=x+\OO{x}$, for $x
The last inequality follows by condition [CLT.Mb].
The same result may be achieved by applying a modified logarithmic
Sobolev inequality to the left tail.

Finally, $\sum_{i}(\mathbb{E}\|f_i(X_i)\|^4)^{1/2}=
\sum_{i} \|a_i-\gamma_n\|^2(m_i^0+3(m_i^0)^2)^{1/2}$,
which is of the order of magnitude of ${O} (d_{\mathcal{H}_n}
)$. This, together with~(\ref{eqlindschwarz}) and~(\ref
{eqsobolevpoisson}) and assumption [CLT.Ma], shows that $U_n$ satisfies
the Lindeberg condition, as stated.

We turn now to consider the conditional limit law. As mentioned above,
Theorem 2.1. in~Morris \cite{MR0370871} holds true also for
multi-dimensional variables. We only need to replace condition $(2.9)$
in~Morris \cite{MR0370871} by a multi-dimensional version.
Specifically, we show that
%
%
%e52 ###
\begin{equation}\label{eqmorriscondition}
\lim_{r\rightarrow0}\sup_n\sup_v
\mathbb{E} \biggl\|\sum_{i}[f_i(L_i+M_i)-f_i(L_i)]\biggr\|^2=0,
\end{equation}
where $L_n = (L_1, \ldots, L_{I_n})$ and $M_n = (M_1, \ldots, M_{I_n})$
are Multinomial random variables with probability vector $\bolds{\pi
}^0_n$, and sample sizes $N_n+v_n N_n^{1/2}$ and $rN_n^{1/2}$,
%$(\pi_n^0,I_n, N_n+v_n N_n^{1/2})$ and $(\pi_n^0,I_n, rN_n^{1/2})$,
respectively, where the parameters $v_n = {O}(1)$ and $r$ are specified
as in Morris \cite{MR0370871}, Lemma 2.2.
% (the first
%element is the vector of probabilities, the second is the number of
%'cells', and
%the third is the 'number of balls'). \ale{Explain better where $v$ and
%$r$ come from and remark that $v$ is bounded!}
Notice that $f_i(L_i+M_i)-f_i(L_i)=(a_i-\gamma_n)M_i$. Thus, %letting
%$M_n = (M_1, \ldots, M_{I_n})$,
%
\[
\biggl\|\sum_{i}(a_i-\gamma_n) M_i \biggr\|^2 = \|\mathrm{A}_n M_n - r
N_n^{-1/2}\mathrm{A}_n \mathbf{m}_n^0 \|
=(M_n -\mathbb{E}M_n)^{\top} B_n (M_n-\mathbb{E}M_n) ,
\]
where $\mathbb{E} M_n = rN_n^{1/2}\bolds{\pi}^0$, and $\mathrm
{B}_n=\mathrm{A}_n^{\top}\mathrm{A}_n$. Taking expectation yields
\begin{eqnarray*}
\mathbb{E} (M_n -\mathbb{E}M_n)^{\top} B_n (M_n-\mathbb{E}M_n) & =&
 r
\sqrt{N_n} \operatorname{tr} \bigl( \mathrm{B}_n \bigl( D_{\bolds{\pi}_n^0} -
\bolds{\pi }_n^0 (\bolds{\pi}_n^0 ) ^{\top} \bigr)\bigr) \\
& = & r \frac{1}{\sqrt{N_n}} \operatorname{tr} \biggl( \mathrm{B}_n \biggl(
D_{\mathbf{m}_n^0} -
\frac{\mathbf{m}_n^0 (\mathbf{m}_n^0 ) ^{\top}}{N_n} \biggr)
\biggr)\\
& = &\mathrm{ O}(1) r \frac{1}{\sqrt{N_n}} ,
\end{eqnarray*}
since
%(\mathbf{m}_n^0 ) ^{\top}}{N_n} ) ) & = & \m{tr}
%( \m{U}_{\mathcal{H}_n} \Sigma^{-1}_{\mathcal{H}_n} \m{U}^{\top}_{
%( D_{\mathbf{m}_n^0} - \frac{\mathbf{m}_n^0 (\mathbf{m}_n^0 ) ^{
%& = & \m{tr} ( \Sigma^{-1}_{\mathcal{H}_n} \m{U}^{\top}_{
%( D_{\mathbf{m}_n^0} - \frac{\mathbf{m}_n^0 (\mathbf{m}_n^0 ) ^{
% \m{U}_{\mathcal{H}_n} )\\
%& = & \m{tr} ( \Sigma^{-1}_{\mathcal{H}_n} \Sigma_{\mathcal{H}_n}
%)\\
%& = & \m{tr}(\m{I}_{\mathcal{H}_n}).
%
\[
\operatorname{tr} \biggl( \mathrm{B}_n \biggl( D_{\mathbf{m}_n^0} - \frac{\mathbf{m}_n^0
(\mathbf{m}_n^0 ) ^{\top}}{N_n} \biggr) \biggr) = \operatorname{tr}
(\mathrm{G}_n\mathrm{G}_n^{\top} ) = {O}(1) .
\]
Therefore,
\[
\mathbb{E} \biggl\|\sum_{i}[f_i(L_i+M_i)-f_i(L_i)]\biggr\|^2
= {O}(1) \frac{r}{\sqrt{N_n}} \rightarrow0 ,
\]
which shows that condition~(\ref{eqmorriscondition}) holds, and the
statement in part 1 is proved.

Part 2 of the theorem follows in a straightforward way from the main
theorem in Bentkus \cite{BEN03} and the fact that $\mathbb{E} \| \mathrm
{F}_n^{-1/2} \mathrm{U}^{\top}_{\mathcal{H}_n} (X_{j_n} - \bolds
{\pi}^0_n) \|
^3$ is of order ${O}( d_{\mathcal{H}_n}^{3/2} )$, by the same
arguments used in the proof of part 1.\vspace*{-3pt}
%= rN_n^{1/2} 0( tr(BD_{p^0})-(p^0)^TBp^0 ) .
%Since the entries of $U_{H_n}$ are bounded by $1$, and $p_i^0\leq
%D_1I_n^{-1}$, we can write:
%tr(BD_{p^0})\leq D I_n \sum_{i,j} |(I_{H_n}^{-1})_{ij} |\leq I_n
%d_{H_n}^{3/2} \lambda_{max}(I_{H_n}^{-1}) \leq
%D\frac{I_n^{5/2}}{N_n} ,
%for some (positive) constant $D$. The second inequality is due to
%the following fact: $\sum_{i,j=1}^d|a_{ij}|\leq d^{3/2}
%To see this, note that
%a_{ij}^2)^{1/2}=d \|A\|_F\leq
%d^{3/2}\|A\|=d^{3/2}\lambda_{max}(A) ,
%where $\|A\|_F$ is the Frobenius norm (which satisfies
%$\|A\|_F\leq d^{1/2}\|A\|$). The term $(p^0)^TBp^0$ is exactly
%$\|\gamma\|^2$, and was shown to be upper bounded by
%$Dd_{H_n}/N_n$, for some positive constant $D$. Multiplying the
%above two upper bounds by $rN_n^{1/2}$ we see that the left hand
%side of~(\ref{eqmorrisexpectation}) is of the order of magnitude
%of $\OO{rN_n^{-1/2}I_n^{5/2}}$. Therefore, under the assumption
%$I_n=\oo{N_n^{1/5}}$, underlying the mild scenario,
%condition~(\ref{eqmorriscondition}) holds and the proof is
%complete.
\end{pf*}

The following lemma is a multivariate analog of Lemma 2.1. in~Morris \cite{MR0370871}.\vspace*{-3pt}
\begin{lemma}\label{lmlemma21}
Let $\mathbf{S}_k=(S_{1k},\mathbf{R}_k)=\sum_{i=1}^k \mathbf
{X}_{ik}$, where
$\mathbf{R}_k=(S_{2k},\ldots, S_{pk})$, $\mathbf{X}_{ik}=(X_{i1k},
\mathbf{Y}_{ik})$,
and $\mathbf{Y}_{ik}=
(X_{i2k}, \ldots, X_{ipk)}$. Suppose that $\{\mathbf{X}_{ik}\}
_{i=1}^k$ are independent
random vectors, with $\mathbb{E}X_{i1k}=0, \mathbb{E}\mathbf
{Y}_{ik}=\mathbf{0}$, and
$\operatorname{Var}(\mathbf{S}_k)=I_p$, the $p\times p$ identity matrix. Suppose
$S_{1k}$ satisfies the
uan condition, i.e., $\max_{1\leq i\leq k}\operatorname{Var} X_{i1k}={o}(1)$ as
$k\rightarrow\infty$, and that
$S_{1k}\stackrel{w}{\longrightarrow} N(0,1)$. Finally, suppose that
$\mathbf{R}_k$ satisfies
the (multi-dimensional) Lindeberg condition, i.e., for all $\varepsilon>0$,
\[
\sum_{i=1}^k\mathbb{E} [ \|\mathbf{Y}_{ik}\|^2 ; \|\mathbf
{Y}_{ik}\|^2>\varepsilon]
={o}(1) \qquad (k\rightarrow\infty) .\vspace*{-3pt}
\]
Then $\mathbf{S}_k\stackrel{w}{\longrightarrow} N_p(\mathbf{0},I_p)$.\vspace*{-3pt}
\end{lemma}
\begin{pf}
As in Morris' proof, $S_{1k}$ satisfies the (one-dimensional) Lindeberg
condition,
i.e.,
\[
\sum_{i=1}^k \mathbb{E}[ X_{i1k}^2 ; X_{i1k}^2> \varepsilon
]={o}(1)
\qquad (k\rightarrow\infty) .\vspace*{-3pt}
\]
Therefore,
\begin{eqnarray*}
\sum_{i=1}^k \mathbb{E} [\|\mathbf{X}_{ik}\|^2 ; \|
\mathbf{X}_{ik}\|^2
\varepsilon]&=  &
\sum_{i=1}^k \mathbb{E} [X_{i1k}^2+\|\mathbf{Y}_{ik}\|^2 ;
X_{i1k}^2+\|\mathbf{Y}_{ik}\|^2
\varepsilon] \\[-2pt]
&  \leq&2 \sum_{i=1}^k \mathbb{E} [ \max\{X_{i1k}^2, \|\mathbf
{Y}_{ik}\|^2\}
; \max\{X_{i1k}^2, \|\mathbf{Y}_{ik}\|^2\}
\varepsilon/2 ] \\[-2pt]
&  \leq& 2 \sum_{i=1}^k \mathbb{E}[ X_{i1k}^2 ; X_{i1k}^2>
\varepsilon/2 ]
+ 2 \sum_{i=1}^k\mathbb{E} [ \|\mathbf{Y}_{ik}\|^2 ; \|
\mathbf{Y}_{ik}\|^2>\varepsilon/2]
={o}(1).\vspace*{-3pt}
\end{eqnarray*}
Thus, $\mathbf{S}_k$ satisfies the (multi-dimensional) Lindeberg
condition and the proof
is complete (see, e.g., Bhattacharya and Rao \cite{MR0436272}, pages 183--184).\vspace*{-3pt}
\end{pf}

%Let $\alpha=\lambda_n/I_n$, and let $S=[-N_n,N_n]$. Clearly,
%$\rho=N_n(1-p_i^0)\geq N_n(1-\max_{i\in\mathcal{I}}p_i^0)$. Thus, by
%the assumption
%$\max_{i\in\mathcal{I}}p_i^0\rightarrow0$, we have, for large enough
%$n$, $\rho\geq CN_n$,
%where $C<1$. We get
%CN_n} ,
%which tends exponentially fast (and thus dominates the $I_n$ elements
%of the summation) to zero.

\section*{Acknowledgements}\vspace*{-3pt}

The authors thank Larry Wasserman for his valuable comments, and one
anonymous reviewer and the associate editor for their suggestions,
which greatly improved the exposition and the readability of the
article. This research was supported in part by NSF Grant EIA-0131884
to the National Institute of Statistical Sciences, by NSF Grant
DMS-06-31589, Army contract DAAD19-02-1-3-0389 and a Health Research
Formula Fund Award granted by the Commonwealth of Pennsylvania's
Department of Health.\vadjust{\goodbreak}

% imsref loaded by akundreckaite, 2011-10-11 15:04:36
%

\printhistory

\end{document}